\documentclass[preprint,1p,fleqn]{elsarticle}
\usepackage{amsfonts,amsmath,amsthm,amssymb}

\usepackage{makeidx}         % allows index generation
\usepackage{graphicx}        % standard LaTeX graphics tool
\usepackage{multicol}        % used for the two-column index
\usepackage[bottom]{footmisc}% places footnotes at page bottom

\usepackage{subfigure}
\def\NN{\hbox{I\kern-.2em\hbox{N}}}
\def\RR{{\mathop{{\rm I}\kern-.2em{\rm R}}\nolimits}}

\def\f{{\bf F}}

\def\n{{\bf n}}
\def\x{{\bf x}}
\def\y{{\bf y}}
\def\lam{{\lambda}}

\def\procO{\texttt{\scalebox{.73}[.85]{procedure 1}}}
\def\procT{\texttt{\scalebox{.73}[.85]{procedure 2}}}

\newcommand{\be}{\begin{equation}}
\newcommand{\ee}{\end{equation}}
\newcommand{\ba}{\begin{eqnarray}}
\newcommand{\ea}{\end{eqnarray}}
\newcommand{\supp}{\mathop{\mathrm{supp}}}

\begin{document}

\begin{frontmatter}

\title{A study on spline quasi--interpolation based quadrature rules for the isogeometric Galerkin BEM}

\author[label1]{Antonella Falini}
\address[label1]{INdAM c/o Department of Mathematics and Computer Science, University of
Florence,\\
Viale Morgagni 67, Firenze, Italy }
\ead{antonella.falini@unisi.it}

\author[label1]{Tadej Kandu\v{c}}
\ead{tadej.kanduc@unifi.it}

\begin{abstract}
Two recently introduced quadrature schemes for weakly singular integrals \cite{CFSS18} are investigated in the context of boundary integral equations arising in the isogeometric formulation of Galerkin Boundary Element Method (BEM). 
In the first scheme, the regular part of the integrand is approximated by a suitable quasi--interpolation spline. In the second scheme the regular part is approximated by a product of two spline functions.\\
The two schemes are tested and compared against other standard and novel methods available in literature to evaluate different types of integrals arising in the Galerkin formulation.
Numerical tests reveal that under reasonable assumptions the second scheme convergences with the optimal order in the Galerkin method, when performing $h$-refinement, even with a small amount of quadrature nodes.
The quadrature schemes are validated also in numerical examples to solve 2D Laplace problems with Dirichlet boundary conditions.
\end{abstract}


\end{frontmatter}

%%%%%%%%%%%%%%%%%%%%%%%%%%%%%%%%%%
\section{Introduction}\label{sec:intro}

% section 1: INTRODUCTION
% section 1: INTRODUCTION
Boundary Element Method (BEM) is a numerical technique %which allows 
to transform the differential problem into an integral one, where the unknowns are defined only on the boundary of the computational domain \cite{costabel1986principles,BEMbook}. 
The main two advantages of the method are the dimension reduction of the problem and the simplicity to treat external problems. As a major drawback, the integral formulation involves Boundary Integral Equations (BIEs), which contain singular kernel functions. Therefore,
robust and precise quadrature formulae are necessary to provide an accurate numerical evaluation. The solution of the considered BIE is then obtained by collocation or Galerkin procedures.
 
The isogeometric formulation of boundary element method (IgA-BEM) has been successfully applied to 2D and 3D problems,
such as linear elasticity \cite{Simp2012}, fracture mechanics \cite{peng2017isogeometric}, acoustic \cite{TauRodHug} and Stokes flows \cite{heltai2014nonsingular}. 
%In all the experimental tests, the isogeometric approach revealed a superior
%accuracy with respect to FEM per used degrees of %freedom.
 
Recently, the IgA paradigm has been combined for the first time to the Symmetric Galerkin Boundary Element Method (IgA-SGBEM) \cite{ADSS1, ADSS2, Nguyen16}, which has revealed to be very effective among BEM schemes. Moreover, the full potential of B-splines over the more common Lagrangian basis has been recently exploited in \cite{ACDS3}.
\medskip

In this work we frame the two quadrature procedures in \cite{CFSS18} in a Galerkin IgA-BEM for the 2D Laplace problem with Dirichlet boundary conditions. 
In particular, the derived quadrature formulae are obtained using a \emph{quasi--interpolation} (\emph{QI}) operator, firstly introduced in \cite{MSbit09} and then applied to construct quadrature rules for regular integrals in \cite{MSJcam12}. The second procedure has been successfully applied in a Galerkin adaptive BEM using hierarchical B-splines in \cite{QIBEM2018}. The authors also provide some theoretical results about the convergence order of the quadrature rule, when $h$-refinement is performed.

In this paper we experimentally test both procedures in \cite{CFSS18} for %the evaluation of 
the regular and weakly singular integrals occurring in the Galerkin formulation. We compare the achieved accuracy with other quadratures available in literature and suitable for the evaluation of the assayed boundary integrals; namely the methods in \cite{ACDS3, Alpert1999, BWR2017, telles1987self}. Moreover, we recall some results about perturbed Galerkin BEM to provide an estimate for the asymptotic accuracy of the quadratures required to obtain the optimal order of convergence.  
\medskip

The structure of the paper is as follows. Section~\ref{Sub:BEMLaplace} presents the integral formulation for exterior and interior problems. Section~\ref{Sub:B-splines} introduces B-spline representations in conjunction with the domain boundary. The Galerkin formulation in the IgA-BEM context is expressed in Section~\ref{Sub:Galerkin}.  
Section~\ref{sec:quad} is devoted to the quadrature rules. In Section \ref{QI-sum} we recall the basics of the adopted quasi--interpolant and the expression of the derived quadrature rules. In Section \ref{sub:comp} some technical aspects related to the quadratures in BEM context are provided. Section~\ref{sec:acc} deals with the accuracy of the considered quadratures for different types of integrals. In Section~\ref{sec:num} three numerical examples are presented using the Galerkin BEM formulation to model 2D Laplace problems, using the proposed QI based quadratures. Finally, some conclusions follow in Section~\ref{sec:conc}. 

%Section~\ref{sec:quad} introduces the quadrature formulas based on spline quasi--interpolation. In Section~\ref{sec:num} the developed  model is applied to an exterior and two interior 2D Laplace problems, all  suited for adaptivity. Finally, Section~\ref{sec:conc} concludes the paper.
%
%
%\bigskip
%
%QI1, QI2, WR, ALPERT, TELES

%%%%%%%%%%%%%%%%%%%%%%%%%%%%%%%%%%
%\section{The problem}\label{sec:prob}
\section{BEM formulation for interior and exterior Laplace problems}\label{sec:bem}
% SECTION 2: SG BEM \section{Boundary integral model problems} \label{sect:BIE}
In the following Section~\ref{Sub:BEMLaplace} we summarise the main features of the BEM formulation for the 2D Laplace problem with Dirichlet boundary conditions and we derive the considered boundary integral equations. Then,
in Section~\ref{Sub:B-splines}, following the isogeometric paradigm, both the boundary representation and the approximate solution are expressed in a B-spline basis. In Section~\ref{Sub:Galerkin}, the governing boundary integrals for the Galerkin discretization are recalled.

\subsection{Direct and indirect BEM formulation for Laplace problems}\label{Sub:BEMLaplace}
In the present work %we analyse 2D Laplace problem with Dirichlet boundary conditions. 
two different domains are considered: open, bounded, simply connected domains $\Omega\subset \RR^2$ and unbounded, external to an open arc $\Gamma$, $\Omega = \RR^2\setminus\Gamma$. The solution of the Laplace problem is to find $u\in H^1(\Omega)$ that satisfies
\begin{equation} \label{Laplace}
\left\{ \begin{array}{rl}
\Delta u=0& {\rm in}\; \Omega\,,\\
u=u_{D}& {\rm on}\; \Gamma\,.
\end{array} \right.
\end{equation}

%\begin{minipage}{0.4\textwidth}
%\centering
%\begin{equation} \label{exterior}
%\left\{ \begin{array}{ll}
%\Delta u=0&{\rm in}\; \RR^2\setminus \Gamma,\\
%u =u_{D}&{\rm on}\; \Gamma\, .
%\end{array} \right. 
%\end{equation}
%\end{minipage}
%\begin{minipage}{0.4\textwidth}
%\centering
%\begin{equation} \label{interior}
%\left\{ \begin{array}{rl}
%\Delta u=0& {\rm in}\; \Omega\,,\\
%u=u_{D}& {\rm on}\; \Gamma\,.
%\end{array} \right.
%\end{equation}
%\end{minipage}\vskip1em
%{\Bl (When you define (1) and (2) you should describe all the new terms in equations...)}

Given the differential problem \eqref{Laplace},
the boundary element method provides an integral formulation, where the unknown is defined only on the boundary of the considered computational domain \cite{costabel1986principles,BEMbook}. In particular, in potential theory we can express the solution $u$ in terms of double layer and single layer potentials using the representation formula (see \cite{costabel1986principles}),

%The boundary element method allows to evaluate the solution $u$ either of \eqref{exterior} or \eqref{interior} at any point $x$ inside the considered domain, by the so called representation formula. 
\begin{equation}\label{rep}
u(\x) = \frac{1}{2\pi}\int_{\Gamma}\partial_{\n_y} U(\x,\y)\, [u(\y)]\, d\Gamma_\y -\frac{1}{2\pi}\int_{\Gamma}U(\x,\y)\, [\partial_{\n_y} u(\y)]\, d\Gamma_\y,\qquad \x\in\Omega.
\end{equation}

The symbol $\partial_\n$ denotes the normal derivative with respect to the exterior unit normal vector $\n$ on $\Gamma$. The harmonic function $u$ is supposed to be regular both in $\Omega$ and in $\RR^2 \setminus \overline{\Omega}$ with different boundary values on both sides of $\Gamma$. The jump of $u$ across $\Gamma$ %, denoted by $[ \cdot ]$, 
is defined as 
\begin{align*}
[u(\x)] := u|_{\RR^2\setminus \Omega} (\x) - u|_{\overline{\Omega}} (\x), \qquad \x\in\Gamma.
\end{align*}
The function ${-\frac{1}{2\pi}U}$ is the fundamental solution of the 2D Laplace operator,
\begin{align*}
U(\x,\y) := \log \|\x-\y\|_2.
\end{align*}
By applying the trace operator to equation \eqref{rep} we can derive a specific \emph{Boundary Integral Equation} (\emph{BIE}) according to certain assumptions considered, when dealing either with an \emph{exterior} or with an \emph{interior} problem.
%By considering either an \emph{exterior} or an \emph{interior} problem, a different trace operator is applied to \eqref{rep} to obtain a specific \emph{Boundary Integral Equation} (\emph{BIE}).
The derived BIE allows us to compute the missing boundary datum. 
%According to which type of problem we are dealing with (i.e., \emph{exterior} or \emph{interior}), applying the trace operator to \eqref{rep} leads to a specific \emph{Boundary Integral Equation} (\emph{BIE}) which allows us to compute the missing boundary datum. 

In the case of an exterior problem, the jump of $u$ across $\Gamma$ vanishes: $[u]\equiv 0$. Therefore, the following BIE can be derived,
\begin{equation}\label{indi}
u_{D}(\x) = -\frac{1}{2\pi}\int_{\Gamma} U(\x,\y)\, \phi(\y)\, d\Gamma_\y, \qquad \x\in\Gamma.
\end{equation} 
In \eqref{indi} the unknown density function  $\phi(\y) := [\partial_{\n_y} u(\y)] \in H^{-1/2}(\Gamma)$ represents the jump of the flux of the solution $u$. 
The space $H^{-1/2}(\Gamma)$ is the dual space of the fractional Sobolev space $H^{1/2}(\Gamma)$, where the duality is defined with respect to the usual $L^2(\Gamma)$-scalar product.
 
In case of an interior problem, we assume $u|_{\RR^2\setminus\Omega}\equiv 0$. The resulting BIE reads
\begin{align}\label{direct}
\frac{1}{2}u_{D}(\x)=\frac{1}{2 \pi} \, \int_{\Gamma}
\partial_{\n_y} U(\x,\y) \, u_{D}(\y)\, d
\Gamma_\y -\frac{1}{2\,\pi}\int_{\Gamma}
U(\x,\y) \, \phi(\y)\, d\Gamma_\y,\quad \x \in \Gamma,
\end{align}
and the unknown function $\phi(\y) :=\partial_{\n_y} u(\y) \in H^{-1/2}(\Gamma)$ is the flux of $u$. 

Thus, we have reformulated the original Laplace problem \eqref{Laplace} in terms of BIE \eqref{indi} and \eqref{direct}.
The formulation in terms of jumps of $u$ is often referred to as \emph{indirect approach}, while the formulation in terms of variables with a clear physical meaning is called the \emph{direct approach}.

Both integrals equations \eqref{indi} and \eqref{direct} are referred to as the Symm's integral equation
\be\label{Symmeq}
V\phi(\x) = f(\x),\quad \x\in\Gamma,
\ee
where $V:H^{-1/2}(\Gamma)\rightarrow H^{1/2}(\Gamma)$ is an elliptic isomorphism and corresponds to the operator
\begin{align*}
V\phi(\x) := - \frac{1}{2 \pi}\int_{\Gamma}U(\x,\y)\phi(\y)\,d\gamma_{\y}.
\end{align*}
The right hand side $f$ in \eqref{Symmeq} is given by $u_D$ in the case of an exterior problem \eqref{indi}, or as 
\begin{align*}
\frac{1}{2}u_{D}(\x)-\frac{1}{2 \pi} \, \int_{\Gamma}
\partial_{\n_y} U(\x,\y) \, u_{D}(\y)\, d
\Gamma_\y
\end{align*}
in the case of an interior problem \eqref{direct}.

\subsection{B-splines}\label{Sub:B-splines}

%{\color{blue} 
%
%Add definition of $S_h$ somewhere.}
%{\color{verde} What space is that one? Where do we need it?}

%The boundary $\Gamma$ is parameterised by a B-spline parametric function $\f : [a,b]\rightarrow \Gamma \subset \RR^2$.

A \emph{knot} vector $ \mathbf T= \left\{t_1,\ldots,t_{N+d +1}\right\}$ is defined as a non-decreasing sequence of knots:
\begin{align*}
%t_i \le t_{i+1}, for i=1,\ldots N+d. and 
t_1 \leq \dots \leq t_{d+1}=:a < t_{d+2} \leq \dots \leq t_N < t_{N+1}=:b \leq \dots \leq t_{N+d+1}.
\end{align*}
The vector $ \mathbf T$ defines a univariate B-spline basis on $[a,b]$ of cardinality $N$ and polynomial degree $d$; the basis is defined by the well-known recursion formula \cite{librodeBoor01}:
 %For the associated partition $\Theta$ on interval $[a,b]\subset\RR$ it holds
%\begin{equation} \label{thetadef}
%a = \theta_1 < \theta_2<\ldots <\theta_L = b,
%\end{equation}
%with $\theta_1=t_{d +1} = a$ and $\theta_L=t_{N+1} = b$. At every breakpoint $\theta_i$, for $i=2,\ldots,L-1$, the corresponding knots are repeated in the inner part of $\mathbf T$ with a multiplicity $m_i$, with $1\le m_i \le d +1$. 
%The $d$ auxiliary knots on the left ($t_1,\ldots,t_{d}$) and on the right ($t_{N+2},\dots,t_{N+d +1}$) may be freely chosen, as long as they preserve the non-decreasing nature of the knot sequence, namely $t_i \le t_{i+1}$, for $i=1,\ldots N+d$. %Note that the following spline dimension formula holds $N= d +1+{\sum_{i=2}^{L-1}}m_i$.
%
%The \emph{B-spline} basis on $\mathbf T$ can be defined with the recursion formula \cite{deBoor01}
\begin{align*}%\label{B_recurrence}
B_{i,0}(t) & := B_{i,0}^{(\mathbf T)}(t) := 
\left\{\begin{array}{ll}
1, & \quad\text{if}\; t_i\le t < t_{i+1},\\
0, & \quad \text{otherwise},
\end{array}\right. \\
B_{i,r}(t) & := \omega_{i,r}(t) B_{i,r-1}(t) + \left(1 - \omega_{i+1,r}(t)\right)B_{i+1,r-1}(t), \quad r=1,\ldots, d,
\end{align*}
where
\begin{align*}
\omega_{i,r} (t) := \left\{
\begin{array}{ll}
\frac{t-t_i}{t_{i+r}-t_i}, & \quad \text{if}\; t_i < t_{i+r},\\
0, & \quad \text{otherwise}.
\end{array}
\right.
\end{align*}
B-splines span a space of splines $\hat S$, whose smoothness depends on the multiplicities of the knots in $\mathbf T$. An interior knot $t_i$ has multiplicity $m\geq1$ if $t_{i-1} < t_i=\dots=t_{i+m-1} < t_{i+m}$ and the space $\hat S$ has a reduced regularity $C^{d-m}$ at the knot $t_i$.

It is common to assume that the boundary $\Gamma$ can be parametrized by a parametric B-spline curve $\f: [a,b] \rightarrow \Gamma \subset \RR^2$, written in the \emph{B-form}, %A planar B-spline curve in $S$ is a parametric curve whose components belong to $S,$ which means it is defined as the image of a vector function $\mathbf{f}: [a,b] \rightarrow \RR^2$,
\begin{align*} %\label{boundarydef}
\f := \sum_{i=1}^N \mathbf{d}_i\, B_{i,d}.
\end{align*}
An ordered set of \emph{control points} in $\RR^2$ is denoted by $\left\{\mathbf{d}_i\right\}_{i=1,\ldots,N}$. 

To recover the interpolation of the first and the last control point, %the left and right auxiliary knots of the extended knot vector are taken coincident with the left and right extrema of the interval
it is common to construct an \emph{open knot vector} by setting $t_1 = \dots = t_{d} = a$ and $t_{N+2} = \dots = t_{N+d+1} = b$. This is the standard choice to define open curves.

For a closed curve geometry, i.e., $\f(a)=\f(b)$, it is more convenient to introduce a periodic definition of the auxiliary knots.
% see, e.g., \cite[Section 10.7]{EssentialsCAGD}. 
To address also the case of a knot vector with multiple knots, we introduce $\rho := d-m+1$, where $m$ is the multiplicity of the knot $a$. In particular, when the knots are simple, $\rho = d$.
The splines are thought to be periodic in a sense that %each B-spline whose support contains $a$ have a counterpart whose support contains $b$
each pair $\{B_{i,d}, B_{N-\rho+i,d}\}$ for $i=1,2,\dots,\rho$
represents one shape function in the physical space. For the periodic compatibility it is sufficient that the $2\rho$ knot differences on the left are identical to the $2\rho$ ones on the right, $t_{i+1} - t_{i} = t_{N-\rho+i+1} - t_{N-\rho+i}$ for $i=1,\dots,2\rho$. Furthermore, $\mathbf d_i = \mathbf d_{N-\rho+i}$ for $i=1,\dots,\rho$. See \cite{EssentialsCAGD}, Section 10.7, for more details.

%{\Bl(}In the present work we analyse recently developed quasi-interpolation based quadratures rules \cite{CFSS18} in the context of Galerkin isogeometric BEM. {\Bl)}
%{\Bl Following the isogeometric paradigm, both the boundary representation and the approximate solution are expressed in a B-spline basis. The boundary $\Gamma$ is parameterised by a B-spline parametric function $\f : [a,b]\rightarrow \Gamma$.}
%%\begin{align*}
%%\f : [a,b]\rightarrow \Gamma.
%%\end{align*}
%{\Bl We restrict our study to $\f \in {\cal C}^2[a,b]$ to avoid problems...
%
%For open geometry we define a B-spline basis on an open knot vector on $[a,b]\subset \RR$...knots...multiplicities $n+1$ for knots $a$ and $b$... Instead, formulation with an extended periodic knot vector is more convenient to describe problems on closed geometries, where $a$ and $b$ map to the same physical point; see, e.g., \cite{EssentialsCAGD}, Section 10.7. %(and Sederberg's CAGD notes (Sec. 6.5) that are not totally correct).
%For that formulation, instead of open knots at the end points, the knots are extended to the left and to the right of the endpoints $a$ and $b$, respectively. Splines are periodic in a sense that each B-spline whose support contains $a$ have a counterpart whose support contains $b$; a pair represents one shape function in the physical space. The periodic compatibility mandates that the $2d-2$ knot intervals closest to $a$ are identical to the $2d-2$ ones at $b$. Furthermore, the first $d$ control points of $\f$ need to coincide with the last $d$ points.}

\subsection{Galerkin formulation} % for IgA-BEM}
\label{Sub:Galerkin}
For both exterior and interior problem, the solution $\phi$ of the considered BIE \eqref{indi} and \eqref{direct} belongs to the Sobolev space $H^{-1/2}(\Gamma)$.
%In the direct and indirect BEM the goal is to find $\phi \in H^{-1/2}(\Gamma)$ that solves BIE \eqref{indi} or \eqref{direct}, respectively. 
The variational formulation of \eqref{Symmeq} is (see \cite{Wendland2}):
\begin{align} \label{weakpb}
\hspace{-0.8cm}\text{\emph{for }} u_D \in H^{1/2}(\Gamma),  \text{\emph{ find }} \phi \in H^{-1/2}(\Gamma)\; \text{\emph{ such that }}\;
{\cal A}(\phi, \psi)={\cal F}(\psi),\; \forall  \psi \in  H^{1/2}(\Gamma),
\end{align}
where the bilinear form ${\cal A}(\phi, \psi)$ and  right-hand side ${\cal F}(\psi)$ are defined as
\begin{align} \label{bilf}
{\cal A}(\phi,\psi):= \int_{\Gamma} \psi({\bf x})\, V\phi({\bf x})\, d\gamma_x,
\qquad
{\cal F}(\psi):= \int_{\Gamma} \psi({\bf x}) \, f({\bf x})\, d\gamma_x.
\end{align}

When using the Galerkin method on \eqref{weakpb}, the infinite dimensional solution space $H^{-1/2}(\Gamma)$ in \eqref{weakpb} is approximated by a finite dimensional subspace $S_h$. 
The parameter $h$ is related to the discretization step size of the subspace $S_{h}$, and the subspace is generated by the lifted B-splines on ${\bf T}$, %Thus  for  exterior problems we set:
\begin{align}\label{eqn:discreteSpace}
S_{h} := 
\left \langle B_{1,d}^{({\bf T})} \circ \f^{-1}, B_{2,d}^{({\bf T})} \circ \f^{-1}, \dots, B_{N,d}^{({\bf T})} \circ \f^{-1} \right \rangle.
\end{align}

%Using Galerkin method, a specific BIE needs to be evaluated according to the chosen approach. In particular we compute an approximation $\phi_h$ of the unknown function $\phi$, expressed in terms of B-spline basis functions of a certain degree $d$. Then, the integration is performed on the support of every basis function in the parameter space.
 
%The unknown $\phi$, either the flux of $u$ or its jump, is expressed in terms of B-spline basis functions. 

Let us introduce coordinates $s,\,t \in [a, b]\subset \RR$ in the parametric domain,
\begin{align*} %\label{st}
s := \f^{-1}(\x), \quad t := \f^{-1}(\y).
\end{align*}
Then the weak form of the exterior problem \eqref{indi} reads
\begin{align} \label{iso-exterior}
& \int_{D_i} B_{i,d}^{(\mathbf T)}(s)\, J(s)\, u_{D}(\f(s))\, ds \nonumber \\
&= -\frac{1}{2 \pi} \int_{D_i} B_{i,d}^{(\mathbf T)}(s)\, J(s)\, \int _{D_i} U(\f(s),\f(t))\, \phi_h(t)\, J(t) \, dt\, ds, \quad i=1,\ldots, N,
\end{align}
where $B_{i,d}^{(\mathbf T)}$ are the test functions of the problem,  $D_i:=\supp B_{i,d}^{(\mathbf T)}$ and $\phi_h\in S_h$ is the approximate solution of $\phi$. The function $J$ denotes the parametric speed of the curve,
\begin{align*}
J(\cdot) := \|\f'(\cdot)\|_2.
\end{align*}
For the interior problem, the corresponding BIE follows from \eqref{direct}, 
\begin{align}\label{iso-interior}
&\frac{1}{2} \int_{D_i} B_{i,d}^{(\mathbf T)}(s)\, J(s)\, u_{D}(\f(s))\, ds \nonumber\\
&= \frac{1}{2 \pi}\int_{D_i} B_{i,d}^{(\mathbf T)}(s)\, J(s)\, \int_{D_i} \partial_{\n_t} U(\f(s),\f(t))\, u_{D}(\f(s))\, J(t)\, dt\, ds \nonumber\\
&-\frac{1}{2 \pi}\int_{D_i} B_{i,d}^{(\mathbf T)}(s)\, J(s)\, \int_{D_i} U(\f(s),\f(t))\, \phi_h(t)\, J(t)\, dt\, ds, \quad i=1,\ldots,N.
\end{align}

To separate the geometrical influence from the singular contribution of the kernel $U(\f(s),\f(t)) =\log \|\f(s)-\f(t) \|_2$, we split $U$ into two functions: $K_1$ regular and $K_2$ weakly singular, defined as
\begin{align*}
K_1(s,t) = \frac 1 2 \log \frac{\|\f(s) - \f(t)\|^2_2}{\delta^2(s,t)},\qquad
K_2(s,t) =  \log \delta(s,t).
\end{align*}
Thus $U(\f(s),\f(t)) =: K_1(s,t)+K_2(s,t)$. To obtain a regular $K_1$, the function $\delta$ needs to be chosen according to the type of domain. Following the idea from \cite{QIBEM2018}, it is defined as: 
\begin{align}
\delta(s,t) := \left\{
\begin{array}{ll}\label{eqn:delta}
|s - t|, & \textrm{if }\Omega = \RR^2\setminus\Gamma, \\
\displaystyle |s-t| \frac{|(s-t)^2 - \gamma^2|}{\gamma^2}, & \textrm{otherwise},%\textrm{if } \Omega\subset\RR^2,
\end{array}
\right.
\end{align}
with $\gamma = b-a$ being the length of the parametric interval.  

%\begin{align*}
%K_1(s,t)&:= \frac{1}{2}\log[R(s,t)]\\
%R(s,t)&:= \frac{||\f(s)-\f(t) ||^2}{(s-t)^2}\\
%\lim_{s\rightarrow t}R(s,t) &= J(t)^2\\
%K_2(s,t) &:= \log|s-t|\\
%\hat K(s,t) &= K_1(s,t) + K_2(s,t)
%\end{align*}
%{\color{verde} some text should be provided between the formulas...}
By writing the approximate solution $\phi_h\in S_h$ in terms of B-splines as
\begin{align*}
{\phi_h = \sum_{j=1}^N \alpha_j\, B_{j,d}^{(\mathbf T)}} \circ \f^{-1}
\end{align*}
and by substituting $K_1$ and $K_2$ in place of $U$, we can rearrange equations \eqref{iso-exterior}--\eqref{iso-interior} in a linear system $A\boldsymbol{\alpha} = \boldsymbol{\beta}$ with unknowns $\boldsymbol \alpha := (\alpha_j)_{j=1}^N$.
In particular, the matrix $A$ consists of entries $A_{i,j} = A_{i,j}^{(1)} + A_{i,j}^{(2)}$ given as
\begin{align}\label{A1}
A_{i,j}^{(\ell)}&:=  -\frac{1}{2 \pi}\int_{D_i} B_{i,d}^{(\mathbf T)}(s)\, J(s)\; ds \int _{D_i} K_\ell(\f(s),\f(t))\,B_{j,d}^{(\mathbf T)}(t) \, J(t) \;dt, %\qquad i,j=1,\dots,N.
%\end{align}
%\begin{align}
%\label{A2}
%a_{i,j}^{(2)}&:=  -\frac{1}{2 \pi}\int_{D_i} B_{i,d}^{(\mathbf T)}(s)\, J(s)\; ds \int_{D_i} K_2(\f(s),\f(t))\,B_{j,d}^{(\mathbf T)}(t) \, J(t) \;dt, \qquad i,j=1,\dots,N.
\end{align}
for $i,j=1,\dots,N$ and $\ell=1,2$.
The entries of the right hand side $\boldsymbol{\beta} := (\beta_i)_{i=1}^N$ are computed according to the type of problem: 
\begin{itemize}
\item For exterior problems, 
\begin{align}\label{bext}
\beta_i &:= \int_{D_i}B_{i,d}^{(\mathbf T)}(s)\, J(s)\,u_{D}(\f(s))\;ds, \qquad i=1,\dots,N.
\end{align}
\item For interior problems, $\boldsymbol{\beta}$ consists of two terms, $\boldsymbol{\beta} = \frac{1}{2}\boldsymbol{\beta}^{(1)}-\frac{1}{2\pi}\boldsymbol{\beta}^{(2)}$. We compute the entries of
$\boldsymbol{\beta}^{(1)}$ by \eqref{bext}, while the elements of $\boldsymbol{\beta}^{(2)}$  as
\begin{align}\label{bint}
\beta_i^{(2)} &:= \int_{D_i} B_{i,d}^{(\mathbf T)}(s)\,J(s)\;ds\,\int_{D_i} \bar{K}(s,t)\,u_{D}(\f(t))\;dt, \qquad i=1,\dots,N.
\end{align}
\end{itemize}
%
%
%The interior problem \eqref{interior} is treated similarly. The matrix $A$ is described again by entries $a_{i,j}^{(1)}+a_{i,j}^{(2)}$ as in Eq.\eqref{Aext}{\Bl (add (8) or have just one number)}.  Nevertheless, it must be noted that kernel $K_1(s,t)$ becomes a singular function when it is evaluated at $s=a$ and $t=b$, or vice versa. Therefore, for closed geometries we propose a new splitting technique based on a different \emph{distant} function. At the end, we can still express the kernel $\hat K$ as a sum of a regular function $\tilde K_1$ and a weakly singular kernel $\tilde K_2$. Further details will be provided in Section \ref{sec:quad}. 
%%
%The right hand side $\b$ consists of two terms. We compute $\b = \b^{(1)}+\b^{(2)}$, with
%$\b^{(1)}$ given by Eq.\eqref{bext} and $\b^{(2)}$ computed as follows, 
%\begin{align}\label{bint}
%b_i^{(2)} &:= \frac{1}{2 \pi}\int_{a}^b B_{i,d}^T(s)\,J(s)\;ds\,\int_{a}^b \bar{K}(s,t)\,u_{D}(\f(t))\;dt, \qquad \forall i.
%\end{align}
Note that the kernel function $\bar{K}(s,t):= \partial_{\n_t} U(\f(s),\f(t))\, J(t)$ is regular everywhere
%$$
%\bar{K}(s,t) := \partial_{\n_t}U(\f(s),\f(t))\,J(t),
%$$
assuming $\f = (F_1,F_2)$ to be ${C}^2$ smooth. More precisely,
\begin{align*}
\lim_{s\rightarrow t}\bar{K}(s,t) = \frac{F_1^{\prime}(t)\, F_2^{\prime\prime}(t)  - F_2^{\prime}(t)\, F_1^{\prime\prime}(t)}{J^{2}(t)},
\end{align*}
see \cite{ACDS3} for further details.

\section{Quadrature rules}\label{sec:quad}

% SECTION 3: QUADRATURES

In this section we summarise the two spline quasi--interpolation based quadrature rules introduced in \cite{CFSS18}. The quadratures are adopted to evaluate regular and singular integrals that appear in the system matrix \eqref{A1} and in the right-hand side vector \eqref{bext}--\eqref{bint}.
Some implementation aspects are explained afterwards.

\subsection{The QI based schemes}\label{QI-sum}
The core idea of the quadrature procedures is to express the integrand functions in terms of simpler functions that can be efficiently integrated. Double integrals in \eqref{A1} are split into two single ones. Regular non-piecewise polynomial parts are approximated by particular quasi--interpolation splines. No special treatment is needed at and near singularities of the kernel $K_2$. Moreover, to simplify and speed up the implementation of the quadratures we can consider the quadrature nodes to be uniformly spaced on the integration domain.

To evaluate entries in \eqref{A1}--\eqref{bint}, the computation of the following two types of integrals needs to be addressed,
\begin{align}\label{reg-int}
I_{B_i}[g] &:= \int_{D_i} B_{i,d}^{(\mathbf T)}(t)\, g(t)\, dt,\\
\label{sing-int}
I_{w_i^s}[g] &:= \int_{D_i} K_2(s,t)\, B_{i,d}^{(\mathbf T)}(t)\, g(t)\, dt,
\end{align}
where we denote $w_{i}^s := K_2(s,\cdot)\, B_{i,d}^{(\mathbf T)}(\cdot)$ and we assume $g\in C(\overline{D}_i)$. Integrals \eqref{sing-int} are considered weakly singular if $s\in D_i$ and nearly singular if $s\notin D_i$ but the distance between $s$ and $D_i$ is sufficiently small. When the distance is sufficiently large, the integral \eqref{sing-int} is regular and it can be considered as a type of \eqref{reg-int}, where $K_2(s,\cdot)$ is hidden inside $g$. 

We recall the basic ideas of the spline QI quadrature procedures in \cite{CFSS18}. In the so-called $\procO$ the whole product $\tilde g := B_{i,d}^{(\mathbf T)}\, g$ is approximated by a quasi--interpolant spline $\sigma_{\tilde g}$ of a chosen degree $p$.
The QI space is defined on an open knot vector ${\boldsymbol\tau}:=\{\tau_{-p},\ldots,\tau_{n+p}\}$, with uniform knots $\tau_{-p}=\ldots=\tau_0<\ldots<\tau_n=\ldots=\tau_{n+p}$ and it is constructed locally on the support $D_i$ of every basis function $B_{i,d}^{(\mathbf T)}$:
\begin{align*}
\hat S_{\boldsymbol\tau}:=\langle B_{-p,p}^{(\boldsymbol\tau)},\ldots,B_{n-1,p}^{(\boldsymbol\tau)}\rangle.
\end{align*}
The breakpoints of $\boldsymbol\tau$ define $n+1$ quadrature nodes localised at $D_i$.
By replacing $\tilde g$ with the quasi--interpolant $\sigma_{\tilde g}$,
\begin{align*}
\sigma_{\tilde g} = \sum_{k=-p}^{n-1} \lambda_{k}(\tilde g)\, B^{({\boldsymbol\tau})}_{k,p} \approx \tilde g,
\end{align*}
where $\lambda_{k}(\tilde g)$ are suitable coefficients, the integrals \eqref{reg-int} and \eqref{sing-int} are approximated by
\begin{align}
\label{eqn:regQ1}
I_{B_i}[g] %\int_{D_i} B_{i,d}^{(\mathbf T)}(t)\, g(t)
&\approx %I_{B_i}^{Q_1}[g]:= %\int_{D_i}  \sum_{k=-p}^{n-1} \lambda_{k}(B_i\, g) B^{\sigma}_{k,p}(t)\, dt = 
\sum_{k=-p}^{n-1} \lambda_{k}(\tilde g)\, \int_{D_i}  B^{(\boldsymbol\tau)}_{k,p}(t)\, dt
= \sum_{k=-p}^{n-1} \lambda_{k}(\tilde g)\, \frac{|\supp B^{(\boldsymbol\tau)}_{k,p}|}{p+1},\\
\label{eqn:singQ1}
I_{w_i^s}[g] 
&\approx %I_{w_i^s}^{Q_1}[g] :=
\sum_{k=-p}^{n-1} \lambda_{k}(\tilde g)\, \int_{D_i}  K_2(s,t)\, B^{(\boldsymbol\tau)}_{k,p}(t)\, dt,
\end{align}
where $|\cdot|$ stands for the size of the region.
In case of integrals \eqref{reg-int} the evaluation gets reduced to the computation of integrals of B-splines \eqref{eqn:regQ1}. For integrals \eqref{sing-int} a preliminary computation of the so-called \emph{modified moments} $\mu_{k,p}^{(i)}(s)$ is needed:
\begin{align}\label{eqn:moments}
\mu_{k,p}^{(i)}(s) := \int_{D_i}  K_2(s,t)\, B^{(\boldsymbol\tau)}_{k,p}(t)\, dt.
\end{align}
Explicit formulae to compute the moments are derived in \cite{ACDS3,CFSS18} for $\delta = |s-t|$. We refer to \cite{QIBEM2018}, when $\delta$ takes the second form in \eqref{eqn:delta}, suitable for a closed boundary curve.

The $\procT$ differs from the $\procO$ in the initial step, where only the function $g$ is approximated by a QI spline $\sigma_g$ in the local space $\hat S_{\boldsymbol\tau}$ of degree $p$:
\begin{align*}
\sigma_{g} = \sum_{k=-p}^{n-1} \lambda_{k}(g)\, B^{(\boldsymbol\tau)}_{k,p} \approx g.
\end{align*}
Thereafter, the product $B_{i,d}^{(\mathbf T)}(\cdot)\, \sigma_g(\cdot)$ %({\color{magenta}$B_{i,d}^{(\mathbf T)}(\cdot)B_{k,p}^{(\boldsymbol\tau)}$ for $k = -p,\ldots,n-1$ \color{blue} I would keep it as I wrote, since explicit products are in the formula below, agree?}) 
of splines is expressed as a linear combination of B-splines basis functions spanning the product space $\Pi$ of degree $d+p$,
\begin{align*}
B_{i,d}^{(\mathbf T)}\, \sigma_{g} = \sum_{k=-p}^{n-1} \lambda_{k}(g)\, B_{i,d}^{(\mathbf T)}\, B^{(\boldsymbol\tau)}_{k,p} = \sum_m \eta_m\, B_{m,d+p}^{(\boldsymbol \tau_\Pi)}.
\end{align*}
The spline space $\Pi$ is defined on a knot vector $\boldsymbol\tau_{\Pi}$ constructed locally on $D_i$ and %$\mathcal G_m$ 
$\eta_m$ are the appropriate coefficients in the new basis. For the details on the construction of the B-spline product space and the representation of the product in the new basis we refer to \cite{Morken91}.
By applying $\procT$, the integrals \eqref{reg-int} and  \eqref{sing-int} are approximated by the following expressions,
\begin{align}
\label{eqn:regQ2}
\hspace{-0.8cm}I_{B_i}[g] &\approx %I_{B_i}^{Q_2}[g] := 
\int_{D_i} \sum_{k=-p}^{n-1} \lambda_{k}(g)\, B_{i,d}^{(\mathbf T)}(t)\, B^{(\boldsymbol\tau)}_{k,p}(t)\, dt
= \sum_m \eta_m\, \frac{|\supp B^{(\boldsymbol\tau_\Pi)}_{m,d+p}|}{d+p+1},\\
\nonumber
\hspace{-0.8cm}I_{w_i^s}[g] %\int_{D_i} K_{2}(s,t)\, B_{i,d}^{(\mathbf T)}(t)\, g(t)\, dt
&\approx %I_{w_i^s}^{Q_2}[g]:= 
\int_{D_i} K_{2}(s,t) \sum_{k=-p}^{n-1} \lambda_{k}(g)\,  B_{i,d}^{(\mathbf T)}(t)\, B^{(\boldsymbol\tau)}_{k,p}(t)\, dt 
%&= \int_{D_i} K_{2}(s,t) \sum_{k=-p}^{n-1} \lambda_{k}(g) \sum_{m} \mathcal{G}_m\, B^{(\boldsymbol \tau_\Pi)}_{m}(t) \,dt \\
\label{eqn:singQ2}
= \sum_{m} \eta_m\, \int_{D_i} K_{2}(s,t)\, B^{(\boldsymbol \tau_\Pi)}_{m,d+p}(t) \,dt \\
&= \sum_{m} \eta_m\, \mu^{(i)}_{m,d+p}(s).
\end{align}

Clearly, choosing a good QI operator is of fundamental importance to obtain accurate quadrature rules. In our study we adopt the Hermite type QI introduced in \cite{MSbit09} and its derivative free variant \cite{MSJcam12}, already framed in a singular integrals context in \cite{CFSS18}. 

Given a function $g$, the quasi--interpolant spline $\sigma_g$, given by the adopted QI operator, is defined on a knot vector $\boldsymbol\tau$ with $n+1$ breakpoints and it can be written in B-form as
\begin{align}\label{eq:genericQI}
\sigma_g = \sum_{j=-p}^{n-1}\lambda_j(g)\, B_{j,p}^{(\boldsymbol\tau)}.
\end{align}

For the scheme in \cite{MSbit09} the coefficients $\lambda_j(g)$ in \eqref{eq:genericQI} are defined as a suitable linear combination of a local subset of values of $g$ and $g^\prime$ at the spline breakpoints. 
For instance, for $p=2$ they can be computed as 
$$
\begin{array}{ll}
\lam_j(g) &=\frac{1}{2}\left( g(\tau_{j+1})+g(\tau_{j+2})\right)-  \frac{\tau_{j+1}-\tau_j}{4}\left( -g'(\tau_{j+1})+g'(\tau_{j+2})\right),\quad j=-1,\ldots,n-2,\\
\lam_{-2}(g) &=\,g(\tau_0), \qquad \lam_{n-1}(g)=g(\tau_n).
\end{array}
$$
This quasi--interpolation scheme is a projector on the considered spline space and has the optimal approximation order $p+1$ for $g\in C^{p+1}([\tau_0, \tau_n])$.

If the variant scheme \cite{MSJcam12} is chosen, then $g'(\cdot)$ in \eqref{eq:genericQI} are computed using a suitable finite difference formula, which approximates the derivative information. 
Also this scheme has optimal approximation order, but it is not a projector. 

As reported in \cite{MSJcam12}, the introduced QI based quadrature formulae for regular integrals are competitive with respect to other QI based schemes (see for instance \cite{Sabit07,Sabit10}), thanks to the usage of the additional derivative information. %In particular such quadrature rules were applied to evaluate regular integrals. %\cite{MSJcam12}. 
The QI variant based quadratures, developed for singular integrals, also exhibit a competitive or superior behaviour when compared to others QI based schemes, see \cite{CFSS18} for more details. 

\subsection{Implementation aspects}\label{sub:comp}

When evaluating the modified moments, numerical instability issues need to be tackled. In this subsection we provide an approach to overcome this problem. A simple observation that can speed up the computation of the modified moments is given afterwards. Some insights regarding additional restrictions related to $\procO$ in the BEM context are stressed at the end.
\medskip

The derived quadrature techniques \eqref{eqn:singQ1} and \eqref{eqn:singQ2} can be applied to evaluate integrals $I_{w_i^s}[g]$ also in case of regular integrals. Nevertheless, we experimentally observed that the computation of the modified moments in \eqref{eqn:moments} exhibits numerical instability as the distance between $s$ and $D_i$ increases; similarly as it was observed for Legendre polynomial based modified moments \cite{AD2002}. Furthermore, the instability increases at a fixed distance, when $D_i$ gets smaller and when the spline degree $p$ in $\mu_{k,p}^{(i)}(s)$ increases. Therefore, when $I_{w_i^s}[g]$ is regular, it is advisable to adopt regular based quadrature rules \eqref{eqn:regQ1} and \eqref{eqn:regQ2}.

When computing the analytical expression of a regular modified moment in \cite{ACDS3,CFSS18,QIBEM2018} in finite arithmetic, a loss of significance in the computation might occur, since the operations contain addends of similar sizes but different signs. Besides using a higher precision arithmetic, the instability effect can be reduced by a tolerance switch: if $D_i$ is far enough from the singularity $s$, the integrand function $K_2$ is a well-behaved function and the corresponding modified moments can be efficiently evaluated with a quadrature rule for regular integrals.
\medskip

Singular integrals \eqref{sing-int} consist of more involved computational steps than the regular integrals \eqref{reg-int}, mainly due to prerequisite computation of the modified moments. The amount of precomputed values can be greatly reduced if we consider, for example, a uniform mesh and that all the shape functions $B_{i,d}^{(\boldsymbol T)}$ are obtained by shifting one instance. In that setting the value $\mu_{k,p}^{(i)}(s)$ in \eqref{eqn:moments} depends only on the relative position of $s$ with respect to the B-spline factor:
\begin{itemize}
\item if
$D_{\tilde i} = D_i + \tilde s$, then $\mu_{k,p}^{(\tilde i)}(\tilde s) = \mu_{k,p}^{(i)}(s)$,
\item if
$B^{(\boldsymbol\tau)}_{k+\tilde k,p}(\cdot +\tilde s) = B^{(\boldsymbol\tau)}_{k,p} (\cdot)$ for the two basis function constructed on $D_i$, then
$\mu_{k+\tilde k,p}^{(i)}(\tilde s) = \mu_{k,p}^{(i)}(s)$.
\end{itemize}
\medskip
 
%{Probably unclear to a reader!: \color{blue}To obtain convergence of the Galerkin solutions $\phi_h$ to the exact $\phi$ while maintaining a fixed amount of quadrature nodes (on the support of test functions)... In that setting, for every constant function $C$ the exactness on integrals $I_{B_i}[C]$ and $I_{w_i^s}[C]$ needs to be satisfied.}

The asymptotic accuracy of a quadrature rule is affected by the accuracy of the approximated function $g$ in $I_{B_i}[g]$ and $I_{w_i^s}[g]$.
The quasi--interpolation operator in \cite{MSbit09} is a projector, i.e., $\sigma_g = g$ for $g = B_{i,d}^{(\boldsymbol T)}$, if $B_{i,d}^{(\boldsymbol T)} \in \hat S_{\boldsymbol \tau}$, see \cite{MSbit09} for details. The projector property implies $p\geq d$ and $\boldsymbol T$ to be a subset of $\boldsymbol \tau$ on $D_i$. Hence, $I_{B_i}[C]$ and $I_{w_i^s}[C]$ can be computed exactly by $\procO$ for any constant $C$, when the projector operator is used. If the rule is exact only for constant functions, we can expect the asymptotic accuracy $O(h)$ for the QI operator, and the overall accuracy $O(h^2)$ for the corresponding quadrature scheme. More precisely, one additional order is obtained due to reduction of the integration domain in the $h$-refinement procedure. 

To obtain the optimal order for $\procO$, QI splines with higher degree should be employed, at least $p\approx 2d$. Furthermore the knot vector $\boldsymbol \tau$ should have multiple knots to satisfy $B_{i,d}^{(\boldsymbol T)} \in \hat S_{\boldsymbol \tau}$ and a suitable generalization of the operator \cite{MSbit09} should be developed.
Another possible limitation of such quasi--interpolation operator is that the construction of the approximant involves also derivatives of $g$, which might not be available for all the considered integrals (e.g., in the outer integrals of $A$ and in $\boldsymbol \beta$).

The derivative free quasi--interpolant variant \cite{MSJcam12} bypasses the latter limitation and it has the same approximation order as the former operator. However, the latter operator is not a projector and hence it is not applicable in $\procO$ to obtain exact values of $I_{B_i}[C]$ and $I_{w_i^s}[C]$. On the other hand, the derivative free variant is a preferable choice for the $\procT$, where the projector is not needed since the quasi--interpolant is constructed only for $g$.

%SECTION 4:
\section{Accuracy of the quadrature rules for the boundary integrals} \label{sec:acc} 

%{\color{red} Comparison of different quadratures for the inner and outer integrals: plots and some comments on convergence (Q1, Q2, Alpert, maybe WR and Teles). Maybe explaining a procedure to automatically increase the amount of nodes.}

%References to the existing results and papers on the convergence of the QI quadrature rules. Then we need to stress the difference in the type of convergence we are considering here (B-splines in the integrals are changing every iteration).
%
%Maybe mentioning something like this (perhaps we write it already when we are deriving the formula for the new QI to motivate our approach):\\

To measure the accuracy of the derived spline quasi--interpolation quadrature schemes $\procO$ and $\procT$, denoted by $\text{QI}_1$ and $\text{QI}_2$, we perform numerical tests for different types of integrals. The tests comprise of regular \eqref{reg-int} and singular integrals \eqref{sing-int}, that appear in the system matrix $A$ and in the right-hand side vector $\boldsymbol \beta$. The quadratures are tested against some other quadrature rules, suitable for the boundary integrals.

The theoretical accuracy of $\text{QI}_1$ and $\text{QI}_2$ with respect to $n$ was studied in \cite{CFSS18}.
In \cite{QIBEM2018} the analysis on the convergence of $\text{QI}_2$ is studied, when $h$-refinement is performed, the amount of quadrature nodes is kept fixed and the regular part of the integrand is sufficiently smooth.
 We recall that the derived convergence order for $\text{QI}_2$ for regular and singular integrals is $O(h^{p+2})$ and $O(h^{p+2}|\log h|)$, respectively.

\subsection{Perturbed system and Strang's lemma}\label{sec:strang}
%
%Exact Galerkin method convergences with ....
%
%Truncation error of quadratures $\to$ perturbed Galerkin method.
%
%To preserve the same order of convergence, the accuracy of a quadrature needs to be properly controlled.
%
%Strang's lemma and needed accuracy on the matrix elements
%\bigskip \bigskip

In this section we summarise results to estimate the needed asymptotic accuracy of the quadrature rules in BEM, when performing $h$-refinement on the discrete spaces. To keep the results concise, we revise only the relevant steps to derive error bounds for the matrix and right-hand side entries.
%for our purposes and 
We refer to \cite{BEMbook} for a detailed analysis. %In conjunction with the following numerical tests in this section we assume the exact solution is sufficiently regular $\phi \in L^2(\Gamma)$.
\medskip

For a sufficiently smooth boundary and a sufficiently regular exact solution $\phi \in H^{d+1}(\Gamma)$ %and the discretization space $S_h$ in \eqref{eqn:discreteSpace} 
there exists a constant  $C>0$ such that
the following error estimate on the approximate solution $\phi_h \in S_h$ holds
\begin{align}\label{eqn:optimConv}
\| \phi - \phi_h \|_{L^2(\Gamma)} \leq C\, h^{d+1} \|\phi\|_{H^{d+1}(\Gamma)},
\end{align}
with $S_h$ defined in \eqref{eqn:discreteSpace}.
Mainly due to the truncation error of the quadrature rules, the Galerkin solution is usually not computed exactly and a notion of a perturbed Galerkin method needs to be introduced. Sufficiently accurate quadrature rules need to be applied at the discretization step size $h$,  so that the computed solution $\tilde \phi_h$ maintains the optimal convergence order,
\begin{align}\label{eqn:optimConv2}
\| \phi - \tilde \phi_h \|_{L^2(\Gamma)} \leq \tilde C\, h^{d+1} \|\phi\|_{H^{d+1}(\Gamma)}.
\end{align}

As a result of the applied quadrature rules, we denote by $\tilde{\mathcal A}_h$ and $\tilde{\mathcal F}_h$ the perturbed functionals of $\mathcal A$ and $\mathcal F$ in \eqref{bilf}. Let us assume the perturbed Galerkin method to be stable, i.e., it satisfies the discrete inf--sup conditions. Then there exists $\gamma>0$ and sufficiently small $h\leq h_0$ such that
\begin{align*}
\gamma &\leq \inf_{\xi_h\in S_h\setminus\{0\}} \sup_{\psi_h\in S_h\setminus\{0\}} \frac{|{\tilde{\mathcal{A}}}_h (\xi_h, \psi_h)|}{\|\xi_h\|_{L^2(\Gamma)}\, \|\psi_h\|_{L^2(\Gamma)}}, \\
\gamma &\leq \inf_{\psi_h\in S_h\setminus\{0\}} \sup_{\xi_h\in S_h\setminus\{0\}} \frac{|{\tilde{\mathcal{A}}}_h (\xi_h, \psi_h)|}{\|\xi_h\|_{L^2(\Gamma)}\, \|\psi_h\|_{L^2(\Gamma)}}.
\end{align*}

The following estimate from Strang's first lemma (see Section~4.2.4 in \cite{BEMbook}) is obtained
\begin{align}
\nonumber \|\phi - \tilde \phi_h\|_{L^2(\Gamma)}
&\leq \|\phi - \phi_h\|_{L^2(\Gamma)} + \|\phi_h - \tilde \phi_h\|_{L^2(\Gamma)} \\
\nonumber &\leq \|\phi - \phi_h\|_{L^2(\Gamma)} + \gamma^{-1} \sup_{\psi_h\in S_h\setminus\{0\}} \frac{|{\tilde{\mathcal{A}}}_h (\phi_h-\tilde \phi_h, \psi_h)|}{\|\psi_h\|_{L^2(\Gamma)}}\\
\nonumber &\leq \|\phi - \phi_h\|_{L^2(\Gamma)} + \gamma^{-1} \bigg( \sup_{\psi_h\in S_h\setminus\{0\}} \frac{|\mathcal{A}(\phi_h, \psi_h) - {\tilde{ \mathcal{A}}}_h(\phi_h, \psi_h)|}{\|\psi_h\|_{L^2(\Gamma)}}\\
\label{eqn:strangEstim}&\phantom{} \hspace{2.87cm}+  \sup_{\psi_h\in S_h\setminus\{0\}} \frac{|\mathcal F(\psi_h) - \tilde {\mathcal F}_h(\psi_h)|}{\|\psi_h\|_{L^2(\Gamma)}} \bigg).
\end{align}
The first term in the right-hand side of estimate \eqref{eqn:strangEstim} is bounded by \eqref{eqn:optimConv}. To satisfy \eqref{eqn:optimConv2} the remaining two consistency error terms in \eqref{eqn:strangEstim} should be sufficiently small. Thus the functionals $\tilde{\mathcal A}_h$ and $\tilde{\mathcal F}_h$ must be sufficiently good approximations of $\mathcal A$ and $\mathcal F$, respectively.

The estimate \eqref{eqn:strangEstim} helps us to bound the needed accuracy of the system matrix $A_h$ and the right-hand side vector $\boldsymbol \beta_h$ at the discretization step size $h$. Let $\tilde A_{h,ij}$ and $\tilde {\beta}_{h,i}$ be the computed entries of $A_{h,ij}$ and $\beta_{h,ij}$ defined in \eqref{A1}, \eqref{bext} and \eqref{bint}. Then
\begin{align*}
A_{h,i,j} &:= \mathcal A \left( B_{j,d}^{(\boldsymbol T_h)} \circ \f^{-1},\, B_{i,d}^{(\boldsymbol T_h)} \circ \f^{-1} \right), &\qquad
\beta_{h,i} &:= \mathcal F \left( B_{i,d}^{(\boldsymbol T_h)} \circ \f^{-1} \right), \\
\tilde A_{h,i,j} &:= \tilde{\mathcal A}_h \left( B_{j,d}^{(\boldsymbol T_h)} \circ \f^{-1},\, B_{i,d}^{(\boldsymbol T_h)} \circ \f^{-1} \right), &\qquad
\tilde {\beta}_{h,i} &:= \tilde{\mathcal F}_h \left( B_{i,d}^{(\boldsymbol T_h)} \circ \f^{-1} \right).
\end{align*}
Let $\phi_h$ and $\psi_h$ be written in B-form in the basis of $S_h$,
\begin{align*}
{\phi_h = \sum_{j=1}^N \alpha_j\, B_{j,d}^{(\mathbf T)}} \circ \f^{-1}, \qquad
{\psi_h = \sum_{i=1}^N \zeta_i\, B_{i,d}^{(\mathbf T)}} \circ \f^{-1}.
\end{align*}
The norm of $\boldsymbol \alpha$ can be bounded by the norm of $\phi_h$ by applying the stability of the B-spline basis in $L^2$ norm, derived from p-norm B-splines estimates in Section 9.3 in \cite{lyche2008spline}, %(for max norm $C_3^{-1}\, \|\boldsymbol q_h\|_\infty \leq \|\phi_h\|_{L^\infty} \leq \|\boldsymbol q_h\|_\infty$),  
\begin{align*}
K^{-1}\, \|\boldsymbol \alpha\|_2 \leq h^{-1/2} \|\phi_h\|_{L^2(\Gamma)} \leq K'\|\boldsymbol \alpha\|_2,
\end{align*}
for $K,K'>0$. A similar estimates holds true for $\boldsymbol \beta$ and $\psi_h$.

Finally, we obtain the estimate on the first consistency error term in \eqref{eqn:strangEstim},
\begin{align*}
\frac{|\mathcal{A}(\phi_h, \psi_h) - {\tilde{ \mathcal{A}}}_h(\phi_h, \psi_h)|}{\|\psi_h\|_{L^2(\Gamma)}}
&= \frac{|\boldsymbol \zeta^T (A_{h} -  \tilde A_{h}) \boldsymbol \alpha|}{\|\psi_h\|_{L^2(\Gamma)}}\\
&\leq \frac{ \|A_{h} -  \tilde A_{h}\|_2\, \|\boldsymbol \zeta\|_2\, \|\boldsymbol \alpha\|_2}{\|\psi_h\|_{L^2(\Gamma)}}\\
&\leq K^2 h^{-1} \| \phi_h \|_{L^2(\Gamma)} \|A_{h} -  \tilde A_{h}\|_2 \\
&\leq K^2 h^{-2}  \| \phi_h \|_{L^2(\Gamma)} \max_{i,j} |A_{h,i,j} -  \tilde A_{h,i,j}|.
%= \sum_{\tau,t} E_{\tau \times t}(\phi_h, \psi_h) \leq \sum E_max \|\phi_h\| \|\psi_h\|
\end{align*}
A similar estimate applies for the second consistency error term,
\begin{align*}
\frac{|\mathcal{F}(\psi_h) - {\tilde{ \mathcal{F}}}_h(\psi_h)|}{\|\psi_h\|_{L^2(\Gamma)}}
&= \frac{|\boldsymbol \zeta^T ({\boldsymbol \beta}_{h} -  \tilde {\boldsymbol \beta}_{h}) |}{\|\psi_h\|_{L^2(\Gamma)}}\\
&\leq \frac{ \| {\boldsymbol \beta}_{h} -  \tilde {\boldsymbol \beta}_{h} \|_2\, \|\boldsymbol \zeta\|_2}{\|\psi_h\|_{L^2(\Gamma)}}\\
&\leq K h^{-1/2}  \| {\boldsymbol \beta}_{h} -  \tilde {\boldsymbol \beta}_{h} \|_2 \\
&\leq K h^{-1} \max_{i} |{\beta}_{h,i} -  \tilde {\beta}_{h,i}|.
%= \sum_{\tau,t} E_{\tau \times t}(\phi_h, \psi_h) \leq \sum E_max \|\phi_h\| \|\psi_h\|
\end{align*}
The derived estimates give us the following bounds on the accuracy of $\tilde A_h$ and $\tilde{\boldsymbol \beta}_h$,
\begin{align}\label{eqn:systemBound}
 \max_{i,j} |A_{h,i,j} -  \tilde A_{h,i,j}| \leq C h^{d+3}, \qquad
\max_{i} |{\beta}_{h,i} -  \tilde {\beta}_{h,i}| \leq C h^{d+2},
\end{align}
which imply the optimal convergence estimate \eqref{eqn:optimConv2} of the perturbed solution $\tilde \phi_h$.

\subsection{Numerical tests}

In this subsection we test $\text{QI}_1$ and $\text{QI}_2$ and compare them with other quadrature rules.
For each of the rules, the tests are repeated twice, with lower and higher amount of quadrature nodes (in the figures we abbreviate the number of quadrature nodes on the support of a B-spline by ``nod.''). 
In all cases we consider uniform meshes with knot vectors ${\bf T}_h$. The spline degree in $S_h$ is fixed to $d=2$ and $p=2$ for the quasi--interpolation spaces. 

The exact integral values are obtained using the integration solver in Wolfram Mathematica.

\subsubsection{Regular integrals $I_{B_i}$}\label{sec:reg}

In the first test we employ the quadratures on a regular integral $I_{B_i}[g]$, defined in \eqref{reg-int}. This type of integrals appear for example in the right hand side $\boldsymbol \beta$. 
Let $g(t) = 3 \sin(\pi(t + 1)) \cos(t + 1)$. In the test we include also a recently developed B-spline weighted quadrature rule (BWR) \cite{BWR2017,ACDS3}, where the integrand B-spline is thought as a \emph{weight} function. The exactness for these rules is imposed on the chosen test spline space or on a refinement of it. 
The weights of the quadrature rule are computed solving a local band system and the quadrature nodes are chosen a priori such that the Schoenberg-Whitney's conditions hold. 

\begin{figure}[t!]
\centering
\subfigure[Lower amount of nodes]{
{\includegraphics[trim = 0.05cm 0cm .95cm 0.5cm, clip = true, height=4cm]{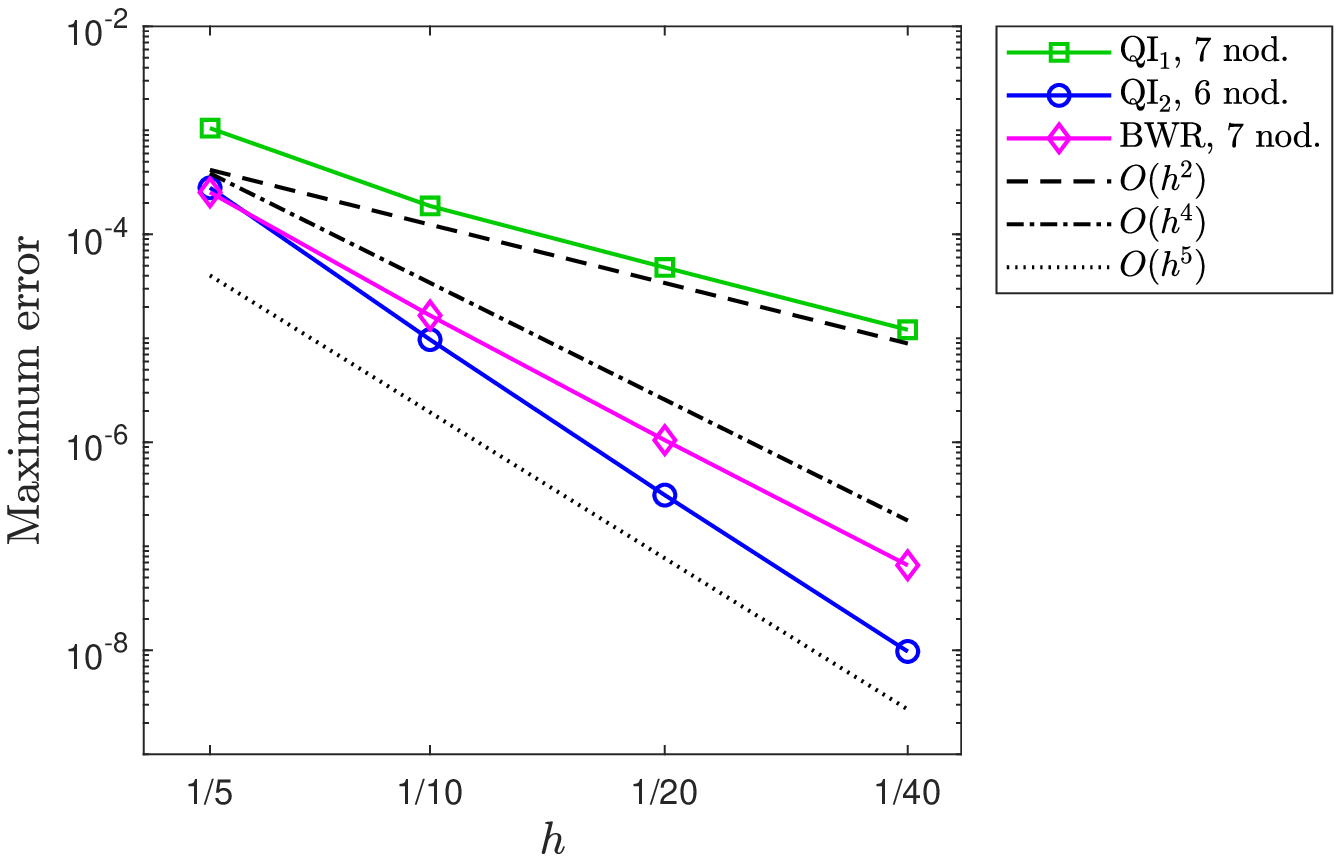}}
}
%\subfigure[Mean error]{
%\includegraphics[trim = 0.5cm 0.5cm 1cm 0.5cm, clip = true, height=6.5cm]{sing_int_conv_mean.eps}
%}
\subfigure[Higher amount of nodes]{
{\includegraphics[trim = 0.05cm 0cm .95cm 0.5cm, clip = true, height=4cm]{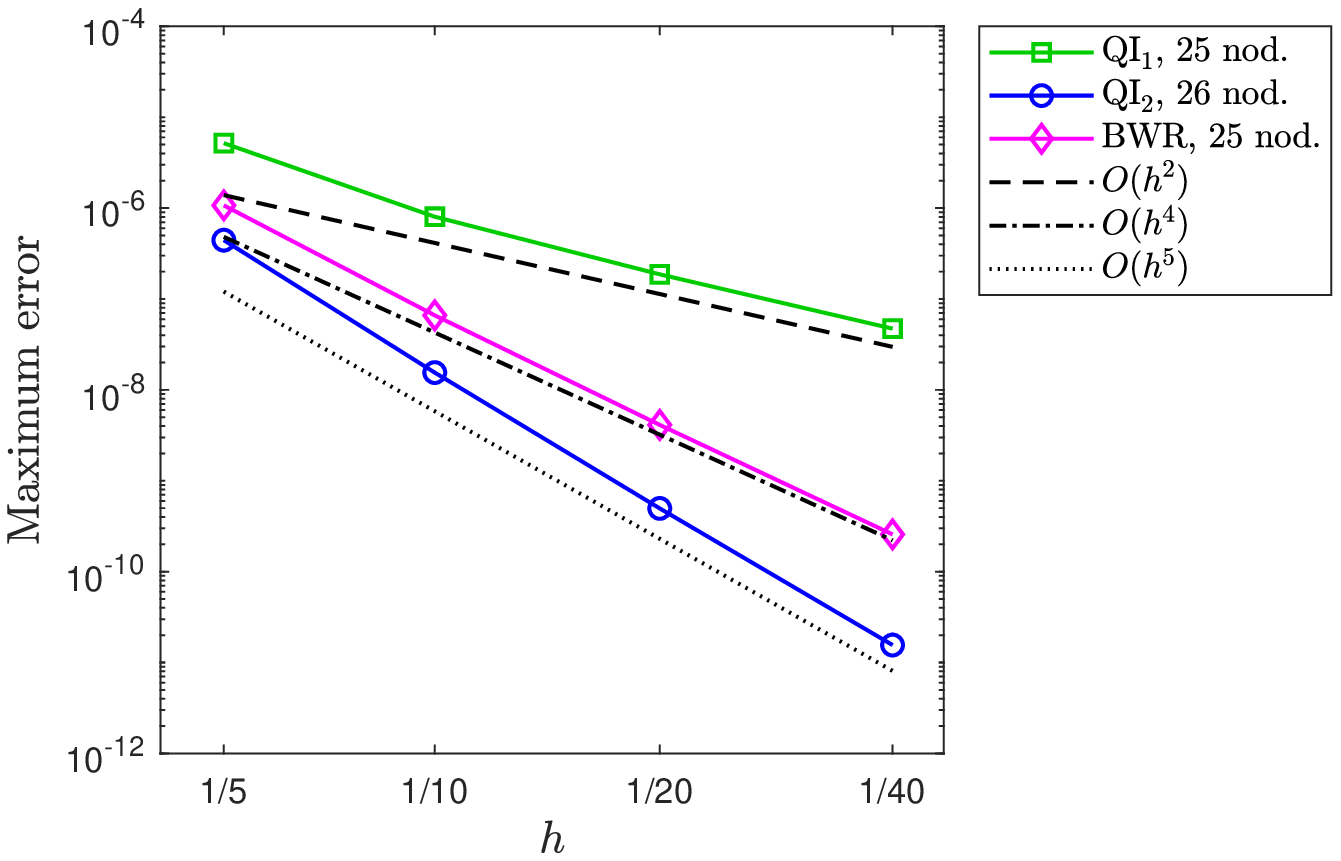}}
}
\caption{Error convergence plots with respect to the mesh size $h$ for a type of integral $I_{B_i}[g]$ in Section~\ref{sec:reg}.}
\label{fig:regIntRHSConvergence}
\end{figure}

%We performed $4$ steps of global refinement on uniformly spaced meshes with $h=1/5,\, 1/10,\, 1/20,\, 1/40$. We computed the maximum error,
For every $h=1/5,\, 1/10,\, 1/20,\, 1/40$ of the uniformly spaced meshes we compute the maximum error of a quadrature scheme,
\begin{align*}
\max_i |I_{B_i}[g]-\tilde I_{B_i}[g]|,
\end{align*}
where $\tilde I_{B_i}[g]$ is the value of a numerically computed integral.
%, where $\star$ indicates either $\text{Q}_1,\text{Q}_2$ or $\text{Q}_{BW}$, when B-spline weighted rules are applied \cite{BWR2017}. 
%The obtained results are shown in 
Figure~\ref{fig:regIntRHSConvergence} reveals the optimal convergence order $O(h^4)$ for BWR, while for $\text{QI}_2$ we observe super convergence $O(h^5)$. Both of the schemes provide the optimal accuracy in the context of the perturbed Galerkin method (see \eqref{eqn:systemBound} in Section~\ref{sec:strang}). The quadrature $\text{QI}_1$ is steadily converging but with a reduced order $O(h^2)$, as discussed in Section~\ref{sub:comp}. As expected, the accuracy of all the rules is improved if we increase the amount of quadrature nodes (Figure \ref{fig:regIntRHSConvergence}(b)). %We used quadratic discretization splines and $p=2$ for both, $\text{QI}_1$ and $\text{QI}_2$.

% Tadej's test for sing integrals' convergence
\subsubsection{Regular and singular integrals $I_{w_i^s}$}\label{sec:inner}

In the second example we measure the accuracy of the integrals with a term $K_2(s,t) = \log |s-t|$. The integrals are regular and singular integrals $I_{w_{h,i}^s}$, with the factor $w_{h.i}^s := K_2(s,\cdot)\, B_{i,d}^{(\mathbf T_h)}(\cdot)$,  that appear as the inner integrals in $A^{(2)}_{i,j}$.
%
%Firstly we analyze the accuracy of our QI {\color{verde}-based} scheme when the following singular integrals (\ref{})
%\begin{align}\label{eqn:singInteg}
%F(h,j,s):=\int_{I} \log|s-t| B_{h,j}^d(t) J(t)\, dt
%\end{align}
%{\color{verde}
%\begin{align}
%\int_{D_i}\log|s-t|B_{j,d}^{({\bf T}_h)}\,J(t)\,dt
%\end{align}}
%are approximated. 
For the test case we consider $g(t)=\sqrt{1+4t^2}$. The quadratic B-splines $B_{i,d}^{(\boldsymbol T_h)}$ are constructed on the interval $[-1,\, 1]$ with uniform open knot vectors ${\bf T}_h$ for the following mesh sizes $h=1/5, 1/10, 1/20, 1/40$. The parameter $s$ is restricted to a priori chosen discrete values. Specifically, $s\in \boldsymbol T_{h/2}$, i.e., $s$ takes the values of all knots in the knot vector ${\bf T}_h$ and all the knot midpoints. 

%{\color{blue}
%Specifically we use 6 knots for the first test and 26 knots for the second one. We stress that the chosen knots are the quadrature nodes used to evaluate \eqref{eqn:singInteg} on $\supp B_{h,j}^d $ for all $j$. }

For every mesh size $h$ we measure the maximum error 
\begin{align*}
%e_{h,{\rm max}} := \max_{j,s} |F(h,j,s) - \tilde F(h,j,s)|,
\max_{i,s} |I_{w_{h,i}^s} - \tilde I_{w_{h,i}^s}|,
\end{align*}
where the computed value of the integral obtained by a quadrature rule is denoted by $\tilde I_{w_{h,i}^s}$. %For QI1 and QI2 {\color{verde} QI-based \texttt{procedure 1} and \texttt{procedure 2} } the parametric speed $J$ in integral \eqref{eqn:singInteg} is approximated by a quadratic QI spline on uniform knots. 

In the test we include also other suitable quadrature rules available in the literature. The hybrid Gauss-trapezoidal quadrature rules, sometimes called just Alpert rules (Alpert), is a class of quadratures for regular and singular functions, that comprise of special quadrature nodes and weights near the (regular or singular) edges of the integration domain \cite{Alpert1999}. The construction exploits a generalization of the Euler-Maclaurin summation formula. A common technique to accurately evaluate weakly singular integrals in BEM is the Telles transformation (Telles), which consists of applying a coordinate transformation to smooth out the singularity and then applying the standard Gaussian quadrature rule to evaluate the regularised integrals \cite{telles1987self}. Finally, we also consider the singular weighted rule (SWR) \cite{ACDS3}. Like $\text{QI}_1$ and $\text{QI}_2$, SWR is based on the precomputed modified moments but the construction of the weights requires to solve a global linear system.

\begin{figure}[t!]
\centering
\subfigure[Lower amount of nodes]{
{\includegraphics[trim = 0.05cm 0cm .95cm 0.5cm, clip = true, height=4cm]{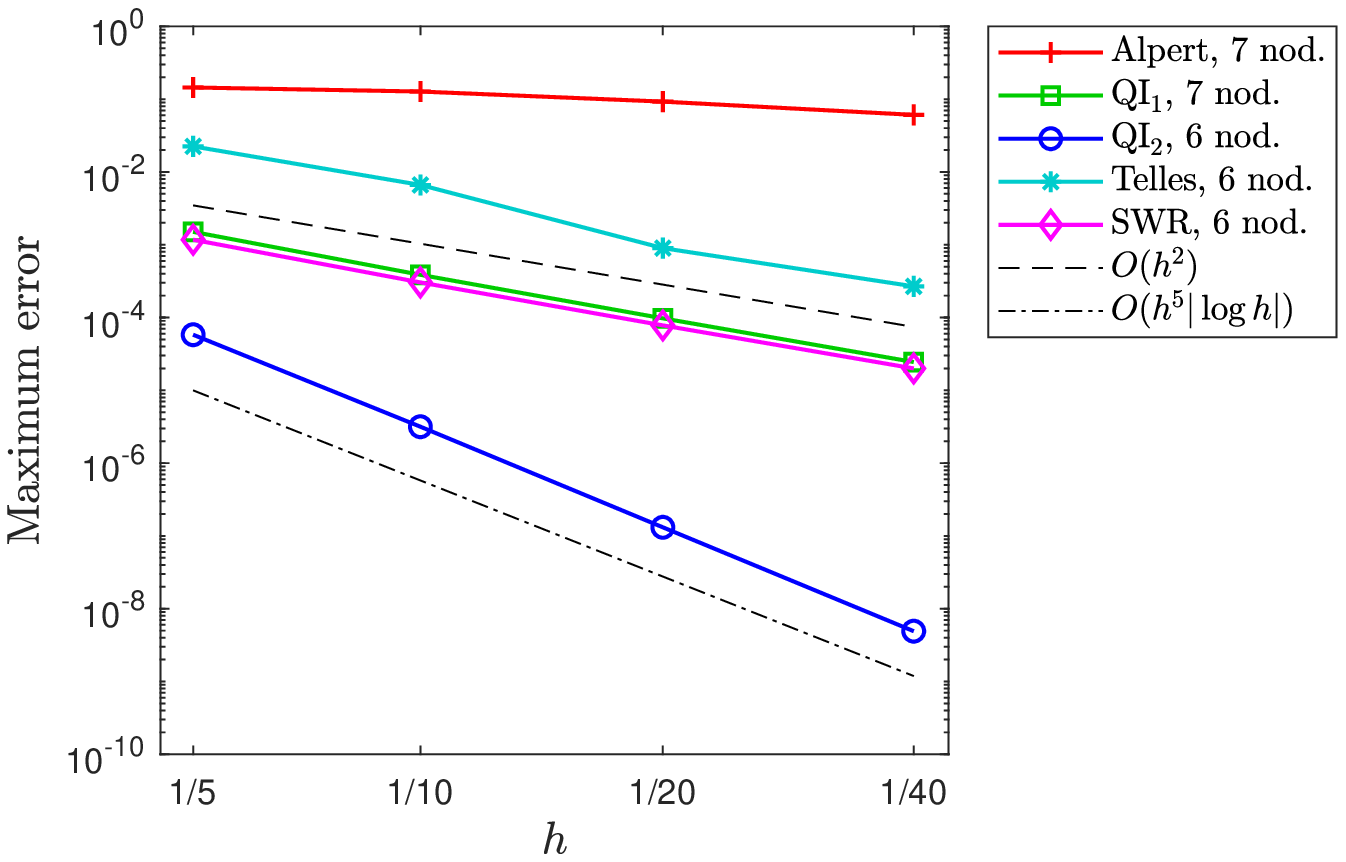}}
}
%\subfigure[Mean error]{
%\includegraphics[trim = 0.5cm 0.5cm 1cm 0.5cm, clip = true, height=6.5cm]{sing_int_conv_mean.eps}
%}
\subfigure[Higher amount of nodes]{
{\includegraphics[trim = 0.05cm 0cm .95cm 0.5cm, clip = true, height=4cm]{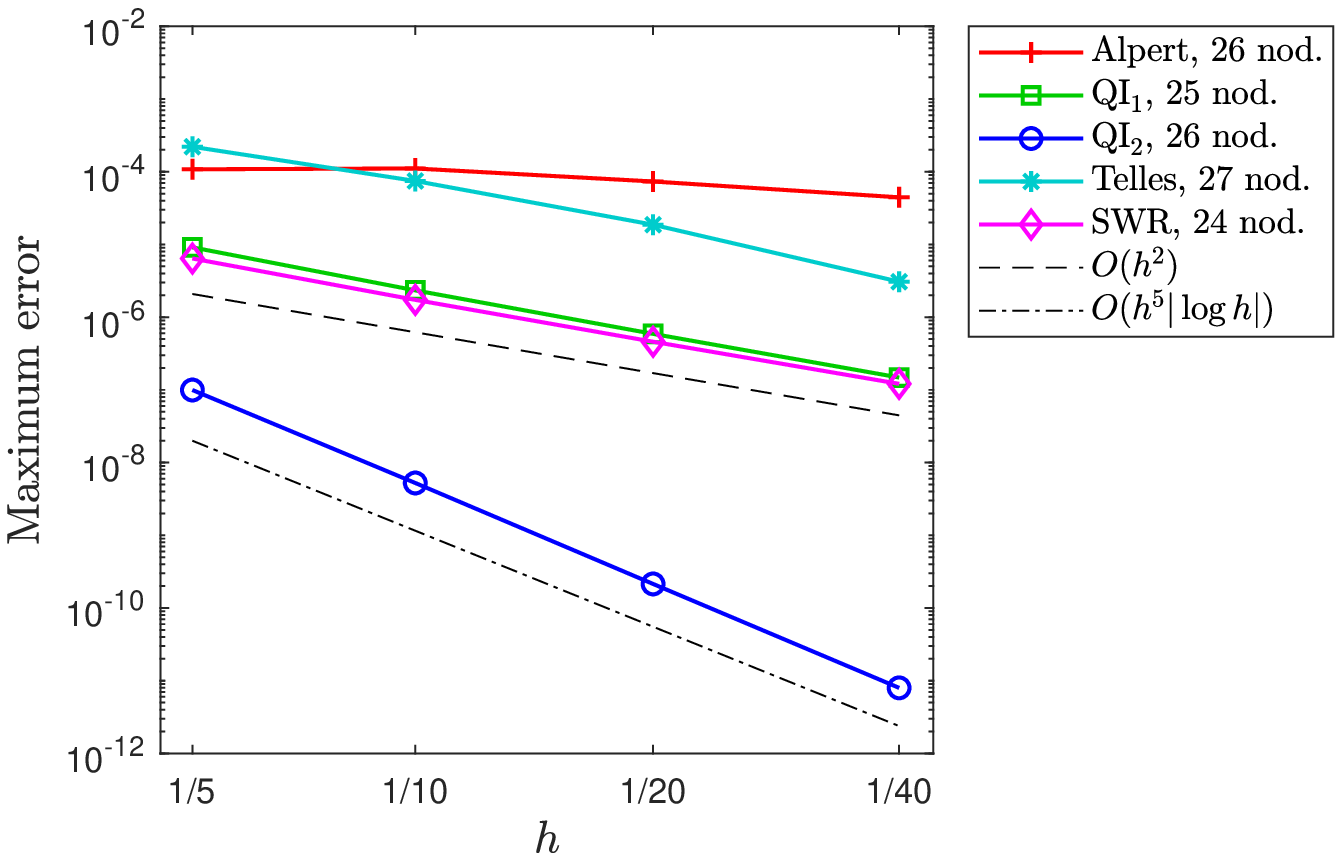}}
}
\caption{Error convergence plots with respect to the mesh size $h$ for the regular and singular inner integrals $I_{w_i^s}$ in Section~\ref{sec:inner} for various quadrature rules available in literature.}
\label{fig:singIntConvergence}
\end{figure}

By comparing the convergence error plots in Figure~\ref{fig:singIntConvergence}(a) and (b) for a fixed $h$ we can observe that all quadrature rules converge with the increased amount of quadrature nodes. Convergence with respect to $h$ reveals the convergence orders $O(h^2)$ and $O(h^5 |\log h|)$ for $\text{QI}_1$ and $\text{QI}_2$ scheme, respectively. The procedure $\text{QI}_2$ outperforms all other quadratures, and it is the only scheme with the optimal convergence order.

% Tadej's test for outer integrals' convergence
\subsubsection{Outer integrals in matrix $A^{(2)}$}\label{sec:outer}

In the last test we focus on regular integrals of the type
\begin{align}\label{eqn:regInteg}
%A_{h}^{(2)}(i,j):=\int_{I} B_{h,j_1}^d(s) J(s)F(h,j_2,s)\, ds,
A_{h,i,j}^{(2)} := -\frac{1}{2 \pi}\int_{D_i} B_{i,d}^{(\mathbf T_h)}(s)\, I_{w_{h,i}^s}\, ds
\end{align}
that appear as the outer integrals in \eqref{A1} for $J\equiv  1$.
The latter choice allows us to exactly evaluate the term $ I_{w_{h,i}^s}$ in \eqref{eqn:regInteg} by computing the corresponding modified moments \eqref{eqn:moments}. Therefore, the error of the numerical integration to compute $A_{h,i,j}^{(2)}$ is contributed solely by the quadrature to compute the outer integral \eqref{eqn:regInteg}.   %Furthermore $A_1\equiv 0$.
%To isolate the error measurement only to the outer integrals, we set $J\equiv 1$; that way the inner integrals are computed accurately (they are equal to the modified moments described in \ref{}) and 
%
On interval $[-1,\, 1]$ we construct uniform meshes with $h=2/5, 1/5, 1/10, 1/20$. 
 
%If we write
%\begin{align*}
%E_h :=  \max_{j_1, j_2} \left|A_{2,h}(j_1,j_2) - \tilde A_{2,h}(j_1,j_2) \right|,
%\end{align*}
%where $\tilde A_{2,h}(j_1,j_2)$ is the computed entry using a quadrature rule. Then by using coordinate transformation we get
%\begin{align*}
%E_{h/2} = \max_{j_1, j_2} A_{2,h/2}(j_1,j_2) &= \max_{j_1, j_2} \int_{I} B_{h/2,j_1}^d(s) J(s) \int_{I} \log|s-t| B_{h/2,j_2}^d(t) J(t)\, dt\, ds\\
%&= \frac{1}{4} \max_{j_1, j_2} \int_{I} B_{h,j_1}^d(s) J(s/2) \int_{I} \log|s-t| B_{h,j_2}^d(t) J(t/2)\, dt\, ds,
%\end{align*}

In this test we include also the previously introduced B-spline weighted quadrature rule (BWR).

\begin{figure}[t!]
\centering
\subfigure[Lower amount of nodes]{
{\includegraphics[trim = 0.05cm 0cm .95cm 0.5cm, clip = true, height=4cm]{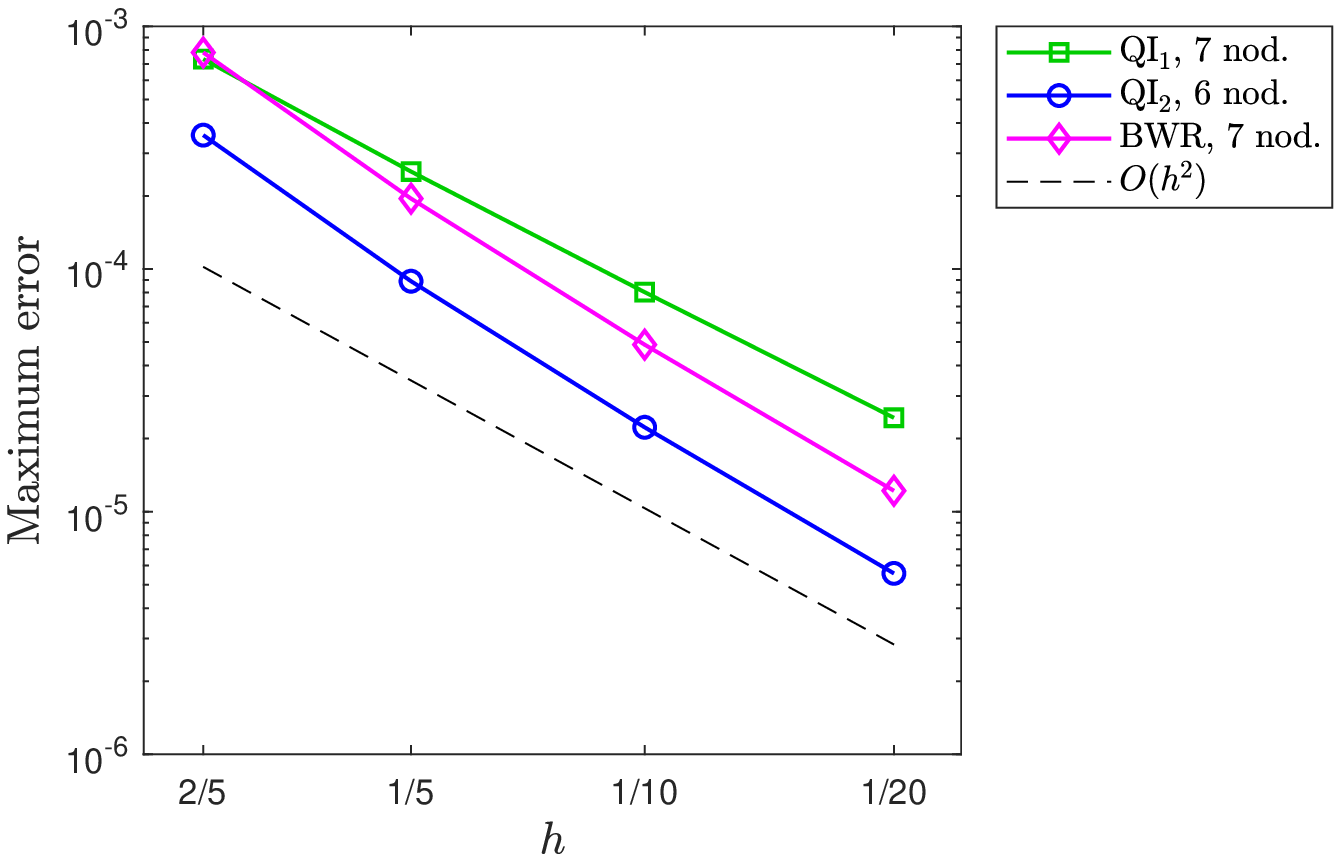}}
}
%\subfigure[Mean error]{
%\includegraphics[trim = 0.5cm 0.5cm 1cm 0.5cm, clip = true, height=6.5cm]{sing_int_conv_mean.eps}
%}
\subfigure[Higher amount of nodes]{
{\includegraphics[trim = 0.05cm 0cm .95cm 0.5cm, clip = true, height=4cm]{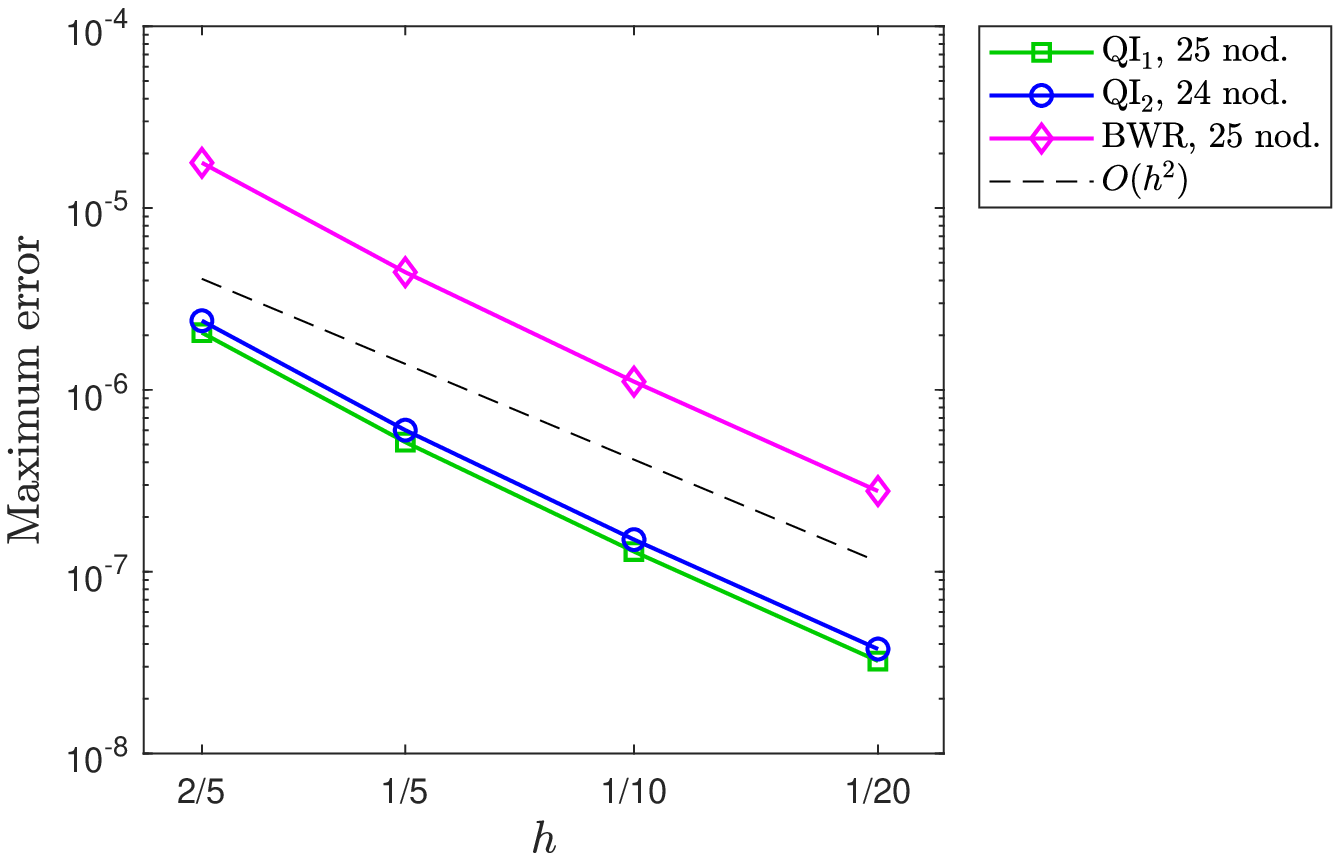}}
}
\caption{Error convergence plots with respect to the mesh size $h$ for the outer integrals in $A^{(2)}$, defined in Section~\ref{sec:outer}.}
\label{fig:regIntConvergence}
\end{figure}

For each mesh size $h$ we measure the maximum error of integrals \eqref{eqn:regInteg} for every $i$ and $j$,
\begin{align*}
%e_{h,{\rm max}}:=\max_{j_1,j_2} \left |A_{2,h}(j_1,j_2) - \tilde A_{2,h}(j_1,j_2) \right|,
\max_{i,j} \left |A^{(2)}_{h,i,j} - \tilde A^{(2)}_{h,i,j} \right|,
\end{align*}
where $\tilde A^{(2)}_{h}(i,j)$ is the computed integral $A^{(2)}_{h}(i,j)$ with a quadrature rule.
Surprisingly, the suboptimal convergence rate is obtained for all the quadratures, even though the function $I_{w_{h,i}^s}$ is well-defined and smooth. From Figure~\ref{fig:regIntConvergence} we can observe the convergence order that is slightly higher than $O(h^2)$. %Note that roughly speaking a quadratic convergence order is obtained just from the reduction of integration area of the integrals when $h$ is reduced. %The obtained order of convergence is suboptimal, since unlike in the first example, the approximated function $F$ varies with $h$. 
Again, the accuracy of the integral approximations is improved, when the amount of quadrature nodes is increased. This can be observed by comparing the plot Figure~\ref{fig:regIntConvergence}(a) for the lower and Figure~\ref{fig:regIntConvergence}(b) for the higher amount of nodes for a fixed $h$.%Hence, refining the space (i.e. considering smaller values for $h$) cannot help to improve the approximation. 

\begin{figure}[t!]
\centering
%\sidecaption[t]
%\subfigure[Lower amount of nodes]{
\includegraphics[trim = 0cm 1cm 1cm 1.5cm, clip = true, height=4cm]{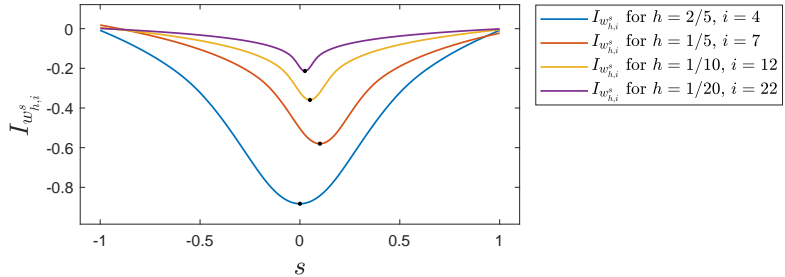}
%}
\caption{Function $I_{w_{h,i}^s}$ for different $h$ and $i$, from the test in Section~\ref{sec:outer}.}
\label{fig:outerIntegFunF}
\end{figure}

To clarify the suboptimal convergence in the test, we plot the function $I_{w_{h,i}^s}$ as a function of $s$ in Figure~\ref{fig:outerIntegFunF} for different fixed $h$ and $i$. For every mesh we plot the function only for the corresponding most central spline element, namely $i=4,7,12,22$. As we can see from the plot, the function $I_{w_{h,i}^s}$ is a smooth function of $s$ but %the complexity of the function does not reduces with smaller $h$. In fact, 
the highest curvature actually increases with smaller $h$ %near its function minima 
(depicted as dots in the figure).
Since the derivatives of $I_{w_{h,i}^s}$ are not bounded when $h\to 0$, the considered quadrature schemes cannot efficiently approximate this type of integrals with a fixed amount of nodes.
%Due to this fact, a suboptimal order of convergence of these kind of integrals is expected for a fixed amount of quadrature nodes. 

%\begin{figure}[t!]
%\subfigure[Maximum error]{
%\includegraphics[trim = 0.5cm 0.5cm 1cm 0.5cm, clip = true, height=6.5cm]{regular_outer_int_conv_max.eps}
%}
%\subfigure[Mean error]{
%\includegraphics[trim = 0.5cm 0.5cm 1cm 0.5cm, clip = true, height=6.5cm]{regular_outer_int_conv_mean.eps}
%}
%\caption{Converge plot of the outer regular integrals with respect to the size $h$.}
%\label{fig:regIntConvergence}
%\end{figure}

%%%%%%%%%%%%%%%%%%%%%%%%%%%%%%%%%%
\section{Numerical simulation with (Galerkin) BEM}
\label{sec:num}
%% SECTION 5: NUMERICAL SECTION

In this section we test the boundary element model to numerically solve three Laplace boundary value problems. For all the examples we evaluate the governing integrals using the presented $\text{QI}_1$ and $\text{QI}_2$ quadrature schemes. In all cases we construct several successive approximate solutions of the problem by performing a dyadic $h$-refinement procedure on uniform meshes. We measure the relative error of an approximated solution against the exact one in $L^2$ norm with respect to degrees of freedom (DoF). The first numerical example is an exterior problem to an open curve, modelled by the indirect formulation. In the next two examples we employ the direct formulation to model interior problems to closed curves.

\subsection{Exterior Dirichlet problem to arc of parabola}\label{ex1}

In this test we focus on the exterior problem described in \cite{ACDS3}, using the BIE \eqref{iso-exterior}. The Dirichlet BVP is defined in the exterior to an arc of parabola $\Gamma \in \RR^2$, parametrized by a quadratic B-spline curve ($d=2$). The transformation map ${\f}(s) = (s, 1 - s^2)$, $s\in[-1,1]$, is described in terms of B-spline basis with the following knot vector $\bf T$ and set of control points $D$,
\begin{align*}
{\bf T} =
(-1, -1, -1, 1, 1, 1),\qquad
D =
\begin{bmatrix}
-1 && 0 && 1\\
0&& 2&& 0
\end{bmatrix}.
\end{align*}

The Dirichlet datum $u_D$ and the exact solution $\phi$ are
\begin{align*}
%(u_D \circ {\f})(s) &= \frac{1}{12 \pi} \left[- (4s^3-9s+7) \log(s^2+2s+2) + (4s^3-9s-7)  \log(s^2-2s+2) \right] \\
%	&+ \frac 1 {6\pi} \bigg[- (4s^3+3s+7) \log(s+1) + (4s^3+3s-7) \log(-s+1) \\
%	&- (12s^2-1) \arctan\left(\frac 2 {s^2}\right) \bigg]
%	+ \frac 2 {9\pi} (12s^2+7),\\
%(\phi \circ {\f})(s) &= \sqrt{1 + 4s^2}.
%
%
&(u_D \circ {\f})(s) \\
&= \frac{- (7-9s+4s^3) \log(2+2s+s^2) - (7+9s-4s^3)  \log(2-2s+s^2)}{12 \pi} + \frac{14+24s^2}{9\pi}\\
&	+ \frac {- (7+3s+4s^3) \log(1+s) - (7-3s-4s^3) \log(1-s) 
	- (-1+12s^2) \arctan\left(\frac 2 {s^2}\right)}{6\pi},\\
&(\phi \circ {\f})(s) = \sqrt{1 + 4s^2}.
\end{align*}

\begin{figure}[t!]
\centering
\subfigure[Dirichlet datum]{
%\fbox{
\includegraphics[trim = 0cm 0cm 0.5cm 0cm, clip = true, height=3cm]{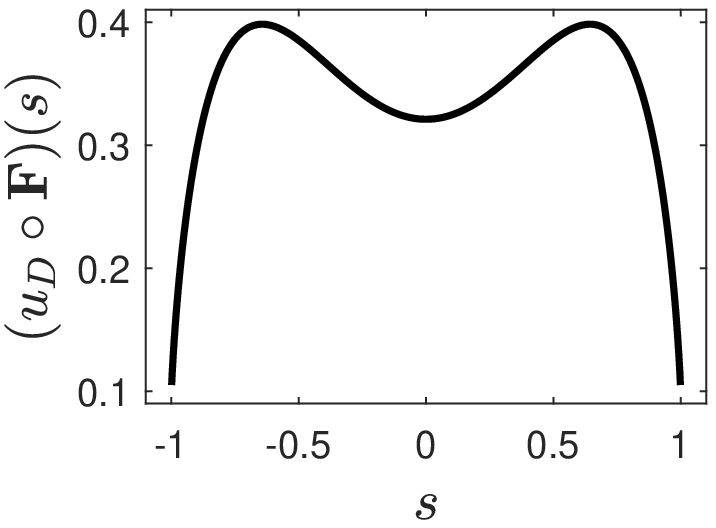}
%}
}
\subfigure[Exact solution]{
%\fbox{
\includegraphics[trim = 0cm 0cm 0.5cm 0cm, clip = true, height=3cm]{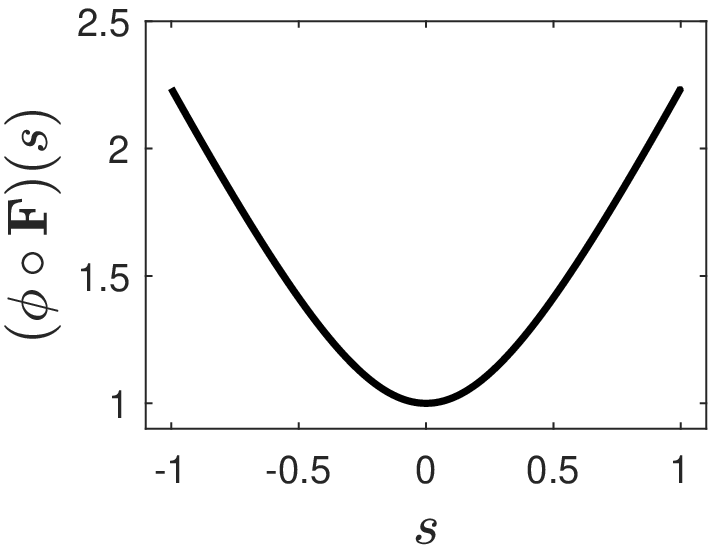}
%}
}
\subfigure[Initial mesh]{
%\fbox{
\includegraphics[trim = 0.25cm 1.75cm 0.5cm 1.75cm, clip = true, height=3cm]{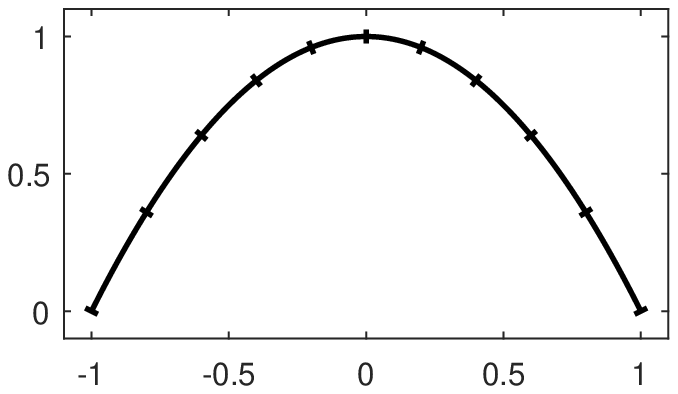}
%}
}
\caption{Parabola test in Section~\ref{ex1}: The Dirichlet datum, the exact solution and the initial mesh. %a) Dirichlet datum. b) Exact solution. c) Initial mesh on the physical geometry.
}
\label{fig:parabola}
\end{figure}
The Dirichlet datum, the exact solution, and the initial mesh in the physical domain are depicted in Figure~\ref{fig:parabola}. Observe that, although the function $u_D \circ \f$ is well defined for $s \in (-1,1)$, its derivative is not bounded when $s\to \pm1$ and can represent an additional limitation for the quadrature $\text{QI}_1$, which needs also the derivative information on the integrand.

The convergence orders of the approximate solutions are shown in Table~\ref{tab:parabola} for the quasi--interpolant degree $p=2$. %with respect to the employed degrees of freedom (DoF). 
Procedure $\text{QI}_1$ has a reduced accuracy near both the edges of the parametric domain. To get the optimal convergence $O(h^{d+1})$ for all refinement steps, a relative high amount of quadrature nodes is needed, %$N_{\rm ref} =
$n+1=49$. On the other hand, it is sufficient to use a small amount of quadrature nodes, $n+1 = 7$, for the procedure $\text{QI}_2$. 
\begin{table}
%\centering
\caption{Parabola test in Section~\ref{ex1}:  Errors and convergence orders of the approximated solutions using quadrature rules $\text{QI}_1$ and $\text{QI}_2$, and $p=2$.} %In both cases we adopt $p=2$ for outer and inner quasi--interpolant splines. The optimal order $O(h^{3})$ is achieved by both procedures, but  with a relative small amount of $7$ nodes for $\text{QI}_2$.}
\begin{tabular}{c||c|c||c|c||}
& \multicolumn{2}{c||}{$\text{QI}_1$, $n+1 = 49$} & \multicolumn{2}{c||}{$\text{QI}_2$, $n+1 = 7$}\\
\hline
DoF & error & conv. & error & conv.\\
\hline
12 & $1.55\cdot 10^{-4}$ &  & $1.57\cdot 10^{-4}$ & \\
22 & $1.66\cdot 10^{-5}$ & 3.69 & $1.66\cdot 10^{-5}$ & 3.70\\
42 & $2.00\cdot 10^{-6}$ & 3.27 & $2.00\cdot 10^{-6}$ & 3.27\\
82 & $2.44\cdot 10^{-7}$ & 3.14 & $2.45\cdot 10^{-7}$ & 3.14\\
162 & $3.06\cdot 10^{-8}$ & 3.05  & $3.05\cdot 10^{-8}$ & 3.06\\
322 & $3.95\cdot 10^{-9}$ & 2.98 & $3.22\cdot 10^{-9}$ & 3.27
\end{tabular}
%\begin{tabular}{l|ccccccc}
%method $\backslash$ DoF & 12 & 22 & 42 & 82 & 162 & 322\\
%\hline
%QI1, 2-2, $N_{\rm ref}=49$& $1.55\cdot 10^{-4}$ &  $1.66\cdot 10^{-5}$ &  $2.00\cdot 10^{-6}$ &  $2.44\cdot 10^{-7}$ &  $3.06\cdot 10^{-8}$ &  $3.95\cdot 10^{-9}$\\
%QI2, 2-2, $N_{\rm ref}=7$& $1.57\cdot 10^{-4}$ &  $1.66\cdot 10^{-5}$ &  $2.00\cdot 10^{-6}$ &  $2.45\cdot 10^{-7}$ &  $3.05\cdot 10^{-8}$ &  $3.22\cdot 10^{-9}$
%\end{tabular}
\label{tab:parabola}
\end{table}

\subsection{Interior Dirichlet problem to a circle}\label{ex2}

To verify the correctness of the model using the direct formulation \eqref{iso-interior}, we consider a simple boundary value problem. The geometry $\Gamma$ is a circle with radius $1/2$, described by the map ${\f}(s) = 1/2\, [\cos(\pi s),\, \sin(\pi s) ]^T$ for $s\in [-1, 1]$. We note that due to the properties of the logarithmic kernel the choice of the radius can greatly influence the condition number of the matrix $A$ and hence the overall results. As analysed and explained in \cite{dijkstra2007relation}, if the radius is equal to one the problem becomes ill-posed and the computed system matrix becomes highly ill-conditioned. 

For the chosen exact solution $u=x$, the Dirichlet datum is $(u_D \circ {\f})(s) = 1/2 \cos(\pi s)$ and the exact flux is $(\phi \circ {\f})(s) = \cos(\pi s)$.

\begin{figure}[t!]
\centering
%\subfigure[Dirichlet datum]{
%\hspace{2cm}
%}
\subfigure[Exact solution]{
%\fbox{
\includegraphics[trim = 0cm 0cm 0cm 0cm, clip = true, height=4.5cm]{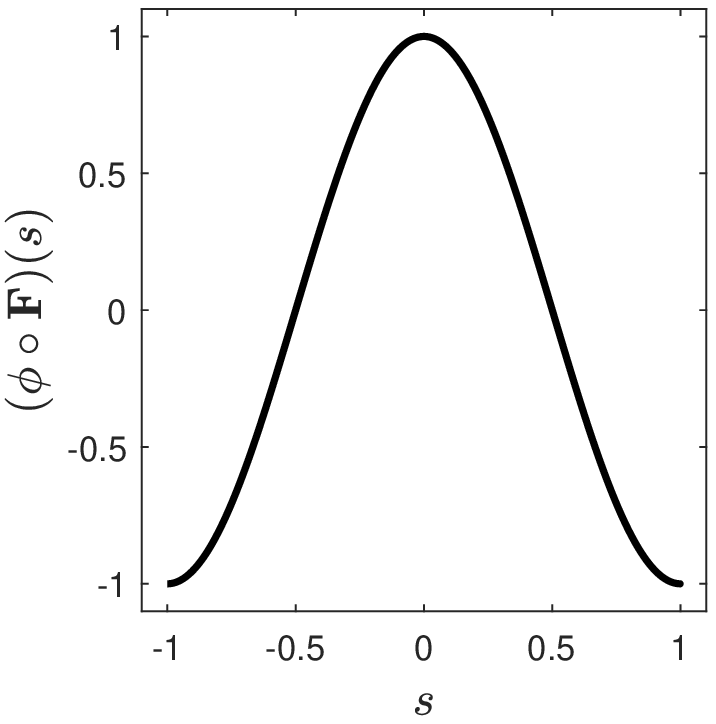}
%}
}
\subfigure[Initial mesh]{
%\fbox{
\includegraphics[trim = 0.cm 0.cm 0.cm 0.cm, clip = true, height=4.5cm]{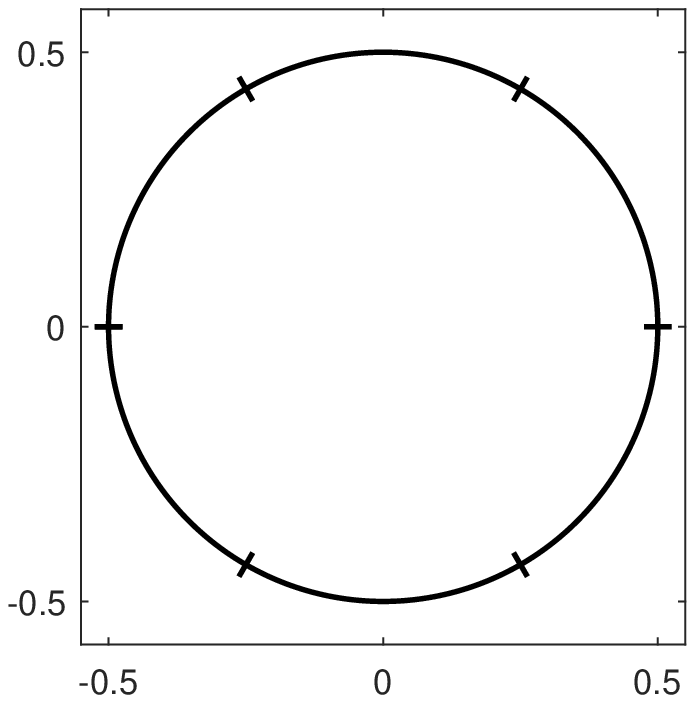}
%}
}

\caption{Circle test in Section~\ref{ex2}: The exact solution and the initial mesh. %a) Dirichlet datum. b) Exact solution. c) Initial mesh on the geometry.}
}
\label{fig:circle}
\end{figure}

The approximate solution is sought in the space of cubic B-splines ($d=3$) with the equally spaced extended knot vector
\begin{align*}
{\bf T} = (-2,\,   -5/3,\,   -4/3,\,  \dots,\, % -1,\,   -2/3,\,   -1/3,\,         0,\,    1/3,\,    2/3,\,    1,\,   
 4/3,\,    5/3,\,    2).
\end{align*}
The exact solution and the initial mesh in the physical domain are depicted in Figure~\ref{fig:circle}.

In Table~\ref{tab:circle} we report errors and convergence rates for the approximate solutions for different amount of DoF. For $\text{QI}_1$ we need to take $p=d$ to satisfy the projector property (see Section~\ref{sub:comp}), whereas for $\text{QI}_2$ it is enough to consider $p=2$.
To obtain the optimal convergence order $4$ for all the refinement steps, we need to considerably increase the amount of quadrature nodes for quadrature $\text{QI}_1$ ($n+1=25$), while for $\text{QI}_2$ we can maintain a small amount of nodes ($n+1 = 5$).

\begin{table}
%\centering
\caption{Circle test in Section~\ref{ex2}: Errors and convergence orders of the approximated solutions using quadrature rules $\text{QI}_1$ for $p=3$ and $\text{QI}_2$ for $p=2$.} %$L^2$ error convergence with respect to DoF. For $\text{QI}_1$ we adopted $p=3$ in outer and inner integrals quasi--interpolant splines, and $n+1=25$ quadrature nodes. For $\text{QI}_2$ we used $p=3$ for the outer quasi--interpolant splines and $p=2$ for the inner ones. Only $n+1=5$ nodes were sufficient to achieve the optimal convergence order $O(h^4)$. {\color{magenta} Why we used $p=3$ for the outer rule in the $\text{QI}_2$?? \color{blue} 3-2 means cubic geometry d=3, p=2.}} 
\begin{tabular}{c||c|c||c|c||}
& \multicolumn{2}{c||}{$\text{QI}_1$, $n+1=25$} & \multicolumn{2}{c||}{$\text{QI}_2$, $n+1=5$}\\
\hline
DoF & error & conv. & error & conv.\\
\hline
6 & $1.66\cdot 10^{-3}$ &  & $1.73\cdot 10^{-3}$ & \\
12 & $7.69\cdot 10^{-5}$ & 4.43 & $8.26\cdot 10^{-5}$ & 4.39\\
24 & $4.40\cdot 10^{-6}$ & 4.13 & $4.67\cdot 10^{-6}$ & 4.14\\
48 & $2.69\cdot 10^{-7}$ & 4.03 & $2.78\cdot 10^{-7}$ & 4.07\\
96 & $1.67\cdot 10^{-8}$ & 4.01 & $1.70\cdot 10^{-8}$ & 4.03\\
192 & $1.05\cdot 10^{-9}$ & 4.00  & $1.05\cdot 10^{-9}$ & 4.02
\end{tabular}
%\begin{tabular}{l|ccccccc}
%method $\backslash$ DoF & 6 &  12 & 24 & 48 & 96 & 192\\
%\hline
%QI1, 3-3, $N_{\rm ref}=25$& $1.66\cdot 10^{-3}$ &  $7.69\cdot 10^{-5}$  &  $4.40\cdot 10^{-6}$ &  $2.69\cdot 10^{-7}$ &  $1.67\cdot 10^{-8}$ &  $1.05\cdot 10^{-9}$\\
%QI2, 3-2, $N_{\rm ref}=5$& $1.73\cdot 10^{-3}$ &  $8.26\cdot 10^{-5}$ &  $4.67\cdot 10^{-6}$ &  $2.78\cdot 10^{-7}$ &  $1.70\cdot 10^{-8}$ &  $1.05\cdot 10^{-9}$\\
%QI2, 3-2, $N_{\rm ref}=9$& $1.66\cdot 10^{-3}$ &  $7.70\cdot 10^{-5}$ &  $4.41\cdot 10^{-6}$ &  $2.69\cdot 10^{-7}$ &  $1.67\cdot 10^{-8}$ &  $1.04\cdot 10^{-9}$
%\end{tabular}
\label{tab:circle}
\end{table}

\subsection{Interior Dirichlet problem to S curve}\label{ex3}

\begin{figure}[t!]
\centering
\subfigure[Exact solution]{
%\fbox{
\includegraphics[trim = 0cm 0cm 0cm 0cm, clip = true, height=4.5cm]{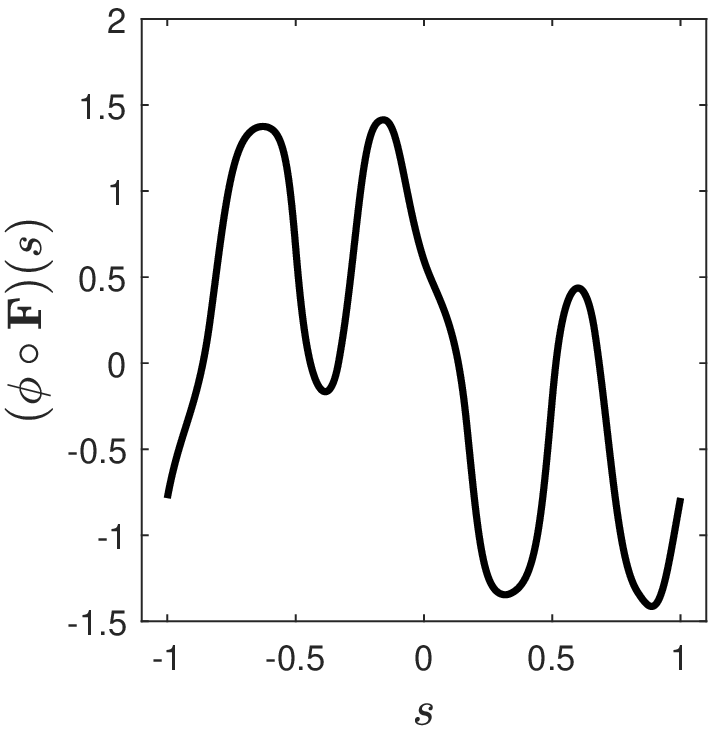}
%}
}
\subfigure[Initial mesh]{
%\fbox{
\includegraphics[trim = 0cm .0cm 0cm 0.0cm, clip = true, height=4.5cm]{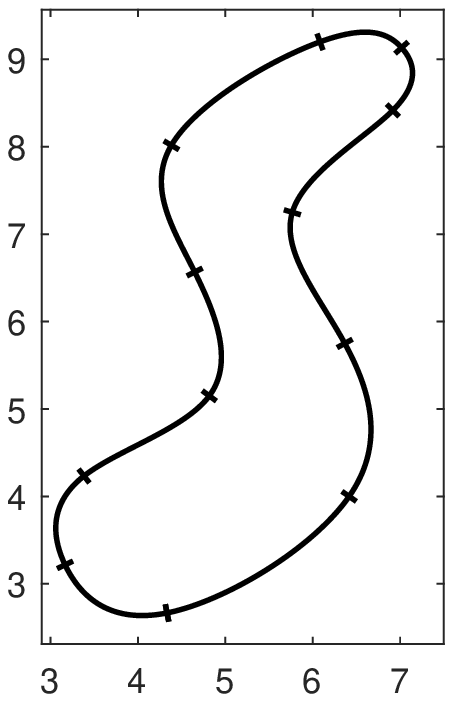}
%}
}
\caption{S curve test in Section~\ref{ex3}: The exact solution and the initial mesh. %a) Initial mesh depicted on the physical domain. b) Exact solution $\phi$.
}
\label{fig:S_curve}
\end{figure}

In the last numerical example we consider a problem with a more involved geometry, a domain described by the closed S curve \cite{ACDS3}. The curve is parametrized by cubic B-splines with the knot vector $\bf T$ and set of control points $D$,
\begin{align*}
{\bf T} &= %(-3/2,\,  -3/2,\,   -3/2,\,   -3/2,\,   -4/3,\,   -7/6,\,   -1,\,   -5/6,\,   -2/3,\,   -1/2,\,   -1/3,\\
%&     -1/6,\,  0,\,  1/6,\,    1/3,\,    1/2,\,    2/3,\,    5/6,\,    1,\,    7/6,\,    4/3,\,    3/2,\,    3/2,\,    3/2,\,    3/2),\\
(-9/6,\,  -9/6,\,   -9/6,\,   -9/6,\,   -8/6,\,   -7/6,\, \dots,\, 7/6,\,    8/6,\,    9/6,\,    9/6,\,    9/6,\,    9/6),\\
D &= 
\begin{bmatrix}
3& 4& 7& 6.5& 5.2& 7.3& 7.1& 6.4& 3.8& 4.7& 5.3& 3& 3& 4& 7\\
3.2& 2.2& 4& 5.8& 7.3& 8.5& 9.2& 9.5& 8& 6.6& 5& 4.3& 3.2& 2.2& 4
\end{bmatrix}.
\end{align*}
The boundary Dirichlet datum is set to $u_D(x,y) = x+y$ and the exact solution reads $(\phi \circ {\f})(s) = (-F'_1(s) + F'_2(s))/\|{\f}'(s)\|_2$. The exact solution and the initial mesh in the physical domain are depicted in Figure~\ref{fig:S_curve}.

In Table~\ref{tab:S_curve} (top) we report errors for the approximate solutions %for different amount of DoF and QI parameters 
when $C^2$ cubic test functions ($d=3$) are used. Again, for $\text{QI}_1$ we set $p=3$ and $p=2$ for $\text{QI}_2$. Since the exact solution $\phi \circ \f$ is only $C^1$ regular, a reduced order of convergence for the approximate solution is expected. This is confirmed by our experiments where the average convergence order is around 2.5.
%\begin{table}
%\begin{tabular}{l|cccccc}
%method $\backslash$ DoF & 12 & 24 & 48 & 96 & 192 & 384\\
%\hline
%QI1, 3-3, $N_{\rm ref}=25$& $1.20{\cdot} 10^{-1}$ &  $3.26{\cdot} 10^{-2}$  &  $4.46{\cdot} 10^{-3}$ &  $6.30{\cdot} 10^{-4}$ &  $1.05{\cdot} 10^{-4}$ &  $1.82{\cdot} 10^{-5}$\\
%QI2, 3-2, $N_{\rm ref}=7$& $1.30{\cdot} 10^{-1}$ &  $3.50{\cdot} 10^{-2}$ &  $4.52{\cdot} 10^{-3}$ &  $6.46{\cdot} 10^{-4}$ &  $1.07{\cdot} 10^{-4}$ &  $1.86{\cdot} 10^{-5}$
%\end{tabular}\\
%
%alternative layout:\\
%\begin{tabular}{l|llllll}
%DoF & 12 & 24 & 48 & 96 & 192 & 384\\
%\hline
%QI1, $N_{\rm r}=25$\\
%error& $1.20{\cdot} 10^{-1}$ &  $3.26{\cdot} 10^{-2}$  &  $4.46{\cdot} 10^{-3}$ &  $6.30{\cdot} 10^{-4}$ &  $1.05 {\cdot} 10^{-4}$ &  $1.82{\cdot} 10^{-5}$\\
%conv.&	&	   1.88&   2.87&   2.82&   2.59&	2.52\\
%\hline
%QI2, $N_{\rm r}=7$\\
%error& $1.30{\cdot}10^{-1}$ &  $3.50{\cdot} 10^{-2}$ &  $4.52{\cdot} 10^{-3}$ &  $6.46{\cdot} 10^{-4}$ &  $1.07{\cdot} 10^{-4}$ &  $1.86{\cdot} 10^{-5}$\\
%conv.&    &	1.89&   2.95&   2.81&   2.59&   2.53
%\end{tabular}
%
%\caption{S curve: $L^2$ error and convergence orders with respect to DoF for $C^2$ cubic test functions.}
%\label{tab:S_curve_cubic}
%\end{table}
%

In Table~\ref{tab:S_curve} (bottom) we can observe an improved convergence order %{\color{magenta} (but then we never reach the optimal order? ... but why...? Why with $C^1$ test functions we don't get 4? boh...) \color{blue} since we use $C^1$ QUADRATICS! It would be interesting to try $C^1$ cubics but i'm not sure if the code can handle it.} 
3 if we employ $C^1$ quadratic test functions ($d=2$) since the test functions have the same regularity as the exact solution. Here $p=2$ for both the quadrature rules. Of course, in this setting the space to describe the geometry is not a subspace of the test space $S_h$.
%\begin{table}
%\begin{tabular}{l|ccccccc}
%method $\backslash$ DoF & 12 & 24 & 48 & 96 & 192 & 384\\
%\hline
%QI1, 2-2, $N_{\rm ref}=13$& $1.24\cdot 10^{-1}$ &  $2.76\cdot 10^{-2}$  &  $2.96\cdot 10^{-3}$ &  $2.55\cdot 10^{-4}$ &  $2.50\cdot 10^{-5}$ &  $2.92\cdot 10^{-6}$\\
%QI2, 2-2, $N_{\rm ref}=7$& $1.26\cdot 10^{-1}$ &  $2.79\cdot 10^{-2}$ &  $2.98\cdot 10^{-3}$ &  $2.58\cdot 10^{-4}$ &  $2.53\cdot 10^{-5}$ &  $2.84\cdot 10^{-6}$
%\end{tabular}
%\caption{S curve: $L^2$ error and convergence orders with respect to DoF for $C^1$ quadratic test functions.}
%\label{tab:S_curve_quadratic}
%\end{table}

To conclude this test, we consider also a case with cubic test functions ($d=3$) that are $C^1$ smooth on the initial mesh by using double knots.  More precisely, the discretization space $S_h$ consists of basis functions that are $C^1$ regular at the initial knots, and $C^2$ continuous at the inserted knots, obtained by dyadic refinements. In this setting, the space to describe the geometry is a subspace of $S_h$ for every step size $h$. 
Note that, in this case for the periodic compatibility $\rho = d-m+1 = 2$ (see Section~\ref{Sub:B-splines}).
The results for the errors and convergence orders are reported in Table~\ref{tab:S_curve2} for $\text{QI}_2$ with $p=2$, $n+1=7$ and $n+1 = 13$. A higher amount of quadrature nodes is necessary to recover the optimal order 4 for the approximate solution. The quadrature scheme $\text{QI}_1$ is not considered in this test since every basis function $B_{i,d}^{(\boldsymbol T)} \in S_h$ should belong to the quasi--interpolation space $\hat S_{\boldsymbol \tau}$ and the involved quasi--interpolation operator cannot handle multiple knots.

\begin{table}
%\centering
\caption{S curve test in Section~\ref{ex3}: Error and convergence orders for $C^2$ cubic (top) and $C^1$ quadratic (bottom) test functions.}
\begin{tabular}{c||c|c||c|c||}
& \multicolumn{2}{c||}{$\text{QI}_1$, $n+1=25$} & \multicolumn{2}{c||}{$\text{QI}_2$, $n+1=7$}\\
\hline
DoF & error & conv. & error & conv.\\
\hline
12 & $1.20\cdot 10^{-1}$ &  & $1.30\cdot 10^{-1}$ & \\
24 & $3.26\cdot 10^{-2}$ & 1.88 & $3.50\cdot 10^{-2}$ & 1.89\\
48 & $4.46\cdot 10^{-3}$ & 2.87 & $4.52\cdot 10^{-3}$ & 2.95\\
96 & $6.30\cdot 10^{-4}$ & 2.82 & $6.46\cdot 10^{-4}$ & 2.81\\
192 & $1.05\cdot 10^{-4}$ & 2.59  & $1.07\cdot 10^{-4}$ &  2.59\\
384 & $1.82\cdot 10^{-5}$ & 2.52 & $1.86\cdot 10^{-5}$ &  2.53
%QI2, quad-quad, 7 nodes& $1.26\cdot 10^{-1}$ &  $2.79\cdot 10^{-2}$ &  $2.98\cdot 10^{-3}$ &  $2.58\cdot 10^{-4}$ &  $2.53\cdot 10^{-5}$ &  $2.84\cdot 10^{-6}$
\end{tabular}
%\hspace{.5cm}
\vspace{0.5cm}

\begin{tabular}{c||c|c||c|c||}
& \multicolumn{2}{c||}{$\text{QI}_1$, $n+1=13$} & \multicolumn{2}{c||}{$\text{QI}_2$, $n+1=7$}\\
\hline
DoF & error & conv. & error & conv.\\
\hline
12 & $1.24\cdot 10^{-1}$ &  & $1.26\cdot 10^{-1}$ & \\
24 & $2.76\cdot 10^{-2}$ & 2.17 & $2.79\cdot 10^{-2}$ & 2.17\\
48 & $2.96\cdot 10^{-3}$ & 3.22 & $2.98\cdot 10^{-3}$ & 3.23\\
96 & $2.55\cdot 10^{-4}$ & 3.53 & $2.58\cdot 10^{-4}$ & 3.53\\
192 & $2.50\cdot 10^{-5}$ & 3.35  & $2.53\cdot 10^{-5}$ &  3.35\\
384 & $2.92\cdot 10^{-6}$ & 3.10 & $2.84\cdot 10^{-6}$ &  3.16
%QI2, quad-quad, 7 nodes& $1.26\cdot 10^{-1}$ &  $2.79\cdot 10^{-2}$ &  $2.98\cdot 10^{-3}$ &  $2.58\cdot 10^{-4}$ &  $2.53\cdot 10^{-5}$ &  $2.84\cdot 10^{-6}$
\end{tabular}
\label{tab:S_curve}
\end{table}

\begin{table}[t]
%\centering
\caption{S curve test in Section~\ref{ex3}: Error and convergence orders for cubic test functions with $C^1$ smoothness at the initial knots, and $C^2$ regularity at the inserted knots.}
\begin{tabular}{c||c|c||c|c||}
& \multicolumn{2}{c||}{$\text{QI}_2$, $n+1=7$} & \multicolumn{2}{c||}{$\text{QI}_2$, $n+1=13$}\\
\hline
DoF & error & conv. & error & conv.\\
\hline
24 & $2.92\cdot 10^{-2}$ &  & $2.65\cdot 10^{-2}$ & \\
36 & $7.12\cdot 10^{-3}$ & 3.48 & $6.67\cdot 10^{-3}$ & 3.40\\
60 & $7.28\cdot 10^{-4}$ & 4.46 & $6.69\cdot 10^{-4}$ & 4.50\\
108 & $4.64\cdot 10^{-5}$ & 4.68 & $3.08\cdot 10^{-5}$ & 5.23\\
204 & $5.51\cdot 10^{-6}$ & 3.35  & $1.73\cdot 10^{-6}$ &  4.53\\
396 & $8.99\cdot 10^{-7}$ & 2.73 & $1.31\cdot 10^{-7}$ &  3.88
%QI2, quad-quad, 7 nodes& $1.26\cdot 10^{-1}$ &  $2.79\cdot 10^{-2}$ &  $2.98\cdot 10^{-3}$ &  $2.58\cdot 10^{-4}$ &  $2.53\cdot 10^{-5}$ &  $2.84\cdot 10^{-6}$
\end{tabular}
\label{tab:S_curve2}
\end{table}

% \begin{table}[H]
% \centering
% \begin{tabular}{|c|l|r|r|r|r|}
% \hline
% & \multicolumn{1}{c|}{Text} & \multicolumn{1}{c|}{Text} & \multicolumn{1}{c|}{Text} & \multicolumn{1}{c|}{Text} & \multicolumn{1}{c|}{text}\\
% \hline
% \parbox[t]{2mm}{\multirow{3}{*}{\rotatebox[origin=c]{90}{rota}}} & text &&&&\\
% & text &&&&\\
% & text &&&&\\
% \hline
% \end{tabular}
% \end{table}

%%%%%%%%%%%%%%%%%%%%%%%%%%%%%%%%%%
\section{Conclusion} \label{sec:conc}

A study of the two recently introduced spline quasi--interpolation quadrature schemes is performed in the context of boundary integral equations in Galerkin IgA-BEM. 
A comparison of the accuracy of the schemes was already done in \cite{CFSS18}, when considering singular integrals. The analysis with respect to the amount of employed quadrature nodes revealed the optimal order of convergence for both approaches. 

In the present paper, numerical tests show a notable difference between the two schemes. For a fixed amount of quadrature nodes the accuracy of the considered integrals is examined, when performing $h$-refinement of the approximation space. The observed rate of convergence is optimal only for the second scheme. 
In the numerical simulations for the 2D Laplace problems, the optimal order of convergence of the approximate solution is achieved with a small number of quadrature nodes, when the second procedure is employed. Regarding the first procedure, the amount of nodes should be increased to recover the optimal order for all the $h$-refinement steps.

In the future work we would like to investigate quadrature schemes for integrals of higher order singularities for more complex differential problems. A valuable contribution would be to derive stable formulae for the modified moments to simplify the construction of the proposed methods.

%
%Numerical tests reveal that under reasonable assumptions the second scheme convergences with the optimal order in the Galerkin method, when performing $h$-refinement, even with a small amount of quadrature nodes.
%The quadrature schemes are validated also in numerical examples to solve 2D Laplace problems with Dirichlet boundary conditions

%Future work: efficient and stable computation of integrals (moments), corners, improving accuracy of the outer rule, multiply-connected domains, mixed problems, 3D.

%%%%%%%%%%%%%%%%%%%%%%%%%%%%%%%%%%
%\section*{Acknowledgements}
%The support by Gruppo Nazionale per il  Calcolo Scientifico
%(GNCS) of the Istituto Nazionale di Alta Matematica (INdAM)
%through ``Progetti di ricerca'' program is gratefully
%acknowledged.

\section*{Acknowledgements}
This work was partially supported by the MIUR “Futuro in Ricerca” programme through the project DREAMS (RBFR13FBI3).
The authors are members of the INdAM Research group GNCS. The INdAM support through GNCS and Finanziamenti Premiali SUNRISE is gratefully acknowledged.

\section*{References}
%\bibliographystyle{plain}
%\bibliography{Dreams_BEM}

\end{document}